\documentclass[leqno, 11pt]{amsart}
\usepackage{amsmath,amssymb,latexsym}
 \newlength{\auxwidth}
 \newlength{\auxheight}
 \newtheorem{definition}{Definition}[section]
 \newtheorem{theorem}[definition]{Theorem}
 \newtheorem{lemma}[definition]{Lemma}
 \newtheorem{proposition}[definition]{Proposition}

 \newtheorem*{theorem*}{Theorem}
\newtheorem*{proposition*}{Proposition}
\newtheorem*{lemma*}{Lemma}

 \theoremstyle{remark}
 \newtheorem{example}[definition]{Example}
 \newtheorem{remark}[definition]{Remark}

  \newtheorem*{acknowledgements}{Acknowledgements}

\newcommand{\op}[1]{\operatorname{#1}}

\newcommand{\acou}[2]{\ensuremath{\langle #1 , #2 \rangle}} 

\newcommand{\tr}{\ensuremath{\op{tr}}}
\newcommand{\Tr}{\ensuremath{\op{Tr}}}
\newcommand{\Tra}{\ensuremath{\op{Trace}}}
\newcommand{\Trace}{\ensuremath{\op{Trace}}}

\newcommand{\Str}{\op{Str}}

\newcommand{\Res}{\ensuremath{\op{Res}}}
\newcommand{\res}{\ensuremath{\op{res}}}

\newcommand{\bint}{\ensuremath{-\hspace{-2,4ex}\int}}

\newcommand{\ind}{\op{ind}}
\newcommand{\Ch}{\op{Ch}}
\newcommand{\Sf}{\op{Sf}}

\newcommand{\End}{\ensuremath{\op{End}}}
\newcommand{\END}{\ensuremath{\op{END}}}
\newcommand{\Cl}{\ensuremath{{\op{Cl}}}}

\newcommand{\C}{\ensuremath{\mathbb{C}}} 
\newcommand{\N}{\ensuremath{\mathbb{N}}} 
\newcommand{\R}{\ensuremath{\mathbb{R}}} 
\newcommand{\Z}{\ensuremath{\mathbb{Z}}}

\newcommand{\UR}{U\times\R}

\newcommand{\Ca}[1]{\ensuremath{\mathcal{#1}}}
\newcommand{\cA}{\Ca{A}}

\newcommand{\cE}{\Ca{E}}
\newcommand{\cH}{\ensuremath{\mathcal{H}}}
\newcommand{\cL}{\ensuremath{\mathcal{L}}}
\newcommand{\cS}{\ensuremath{\mathcal{S}}}
\newcommand{\cD}{\ensuremath{\mathcal{D}}}

\newcommand{\pdo}{$\Psi$DO}
\newcommand{\pvdo}{\ensuremath{\Psi_{\op{v}}}} 
\newcommand{\psido}{$\Psi$DO} 
\newcommand{\psidos}{$\Psi$DO's}

\newcommand{\ev}{{\text{ev}}}
\newcommand{\odd}{{\text{odd}}}

\newcommand{\hotimes}{\hat\otimes}

\newcommand{\sD}{\ensuremath{{/\!\!\!\!D}}}

\newcommand{\sP}{\ensuremath{{/\!\!\!\!P}}}
\newcommand{\sS}{\ensuremath{{/\!\!\!\!\!\;S}}}

\newcommand{\sSE}{\ensuremath{\sS\otimes\cE}}

\begin{document}
\title{A NEW SHORT PROOF OF THE LOCAL INDEX FORMULA\\ AND SOME OF ITS APPLICATIONS} 

\author{Rapha\"el Ponge}

\address{Department of Mathematics, Ohio State University, Columbus, USA.}
\email{ponge@math.ohio-state.edu}

\numberwithin{equation}{section}

\date{}

\begin{abstract}
We  give a new short proof of the index formula of Atiyah and Singer based on combining Getzler's 
rescaling with Greiner's approach of the heat kernel asymptotics. As application we 
can easily compute the Connes-Moscovici cyclic cocycle of even and odd Dirac spectral triples, and then recover the 
Atiyah-Singer index formula (even case) and the Atiyah-Patodi-Singer spectral flow formula (odd case). 
\end{abstract}

\maketitle

The Atiyah-Singer index Theorem~(\cite{AS:IEO1},~\cite{AS:IEO3}) gives a cohomological interpretation of the Fredholm index of 
an elliptic operator, but it reaches its true geometric content in the case Dirac operator for which the index is given
by a local geometric formula. The local formula is somehow as important as the index theorem since, on the one 
hand, all the common geometric operators are locally Dirac operators (\cite{ABP:OHEIT}, \cite{BGV:HKDO}, \cite{LM:SG}, 
\cite{Roe:EOTAM}) and, on the other hand, the local index formula 
is equivalent to the full index theorem (\cite{ABP:OHEIT}, \cite{LM:SG}). 
It was then attempted to bypass the index theorem to prove the local 
index formula. The first direct proofs were made by Patodi, Gilkey,  Atiyah-Bott-Patodi partly by using invariant 
theory (see \cite{ABP:OHEIT}, \cite{Gi:ITHEASIT}). Some years later Getzler (\cite{Ge:POSASIT}, \cite{Ge:SPLASIT}) and 
  Bismut~\cite{Bi:ASITPA}  gave purely analytic proofs, which led to many 
  generalizations of the local index formula (see also~\cite{BGV:HKDO}, \cite{Roe:EOTAM}).

The short proof of  Getzler~\cite{Ge:SPLASIT} combines the Feynman-Kac representation of the heat kernel with an ingenious trick, 
the Getzler rescaling. In this paper we a give a new short proof of the local index 
formula for Dirac operators by combining Getzler rescaling with the 
(fairly standard) Greiner's approach of the heat kernel asymptotics~(\cite{Gr:AEHE}, \cite{BGS:HECRM}). Our proof is 
quite close to other proofs like those by Melrose~\cite[pp.~295-327]{Me:APSIT}, Simon~\cite[Chap.~12]{CFKS:SOAQM}
and Taylor~\cite[Chap.~10]{Ta:PDE2}, but 
here the justification of the convergence of the supertrace of the heat kernel, which is the key of the proof, follows 
from very elementary consideration on Getzler's orders (Lemma~\ref{lem:AS.approximation-asymptotic-kernel}).

In fact, the proof  yields a more general result, for it implies a differentiable version of the asymptotics for the 
supertrace of the heat kernel,  which is hardly accessible by means of a probabilistic representation of the heat kernel as 
in~\cite{Ge:SPLASIT} (see Proposition~\ref{prop:AS.differantiated-local-index-theorem}).
 
 In the second part of the paper we show how this enables us to compute the CM cyclic cocycle~\cite{CM:LIFNCG} 
 associated to a Dirac spectral triple, both in the even case
(Theorem~\ref{thm:even-chern}) and in the odd case (Theorem~\ref{thm:odd-chern}).  Therefore we can bypass the use of 
Getzler's asymptotic pseudodifferential calculus~\cite{Ge:POSASIT} of 
 the previous approaches of the computation of the CM cocycle for Dirac spectral triples 
 (\cite[Remark~II.1]{CM:LIFNCG}; see also \cite{CH:ECCIDO}, \cite{Le:TSPSCI}).    

Recall that the CM cocycle is important because it represents the cyclic cohomology Chern character of a 
spectral triple (i.e. a "noncommutative manifold'') and is given by a formula which is local in the sense of noncommutative geometry 
 (\cite{CM:LIFNCG}; see also Section~\ref{sec.cyclic-chern}). Thus it allows 
 the local index formula to hold in a purely operator theoretic setting. 
For instance, the computation for Dirac spectral triples allows us  to recover, in the even case, the local index formula of Atiyah-Singer  
 and, in the odd case, the spectral flow formula of Atiyah-Patodi-Singer~\cite{APS:SARG3} 
 (\emph{cf.}~\cite[Remark II.1]{CM:LIFNCG} and sections~\ref{sec.even-chern} and~\ref{sec.odd-chern}; 
 see also~\cite{Ge:OCCCHSF} for the odd case).  

The paper is organized as follows. In the first section we recall Greiner's approach of the heat kernel asymptotics 
following~\cite{Gr:AEHE} and~\cite{BGS:HECRM}.  In Section~\ref{sec.getzler} we prove  
the local index formula of Atiyah-Singer and in Section~\ref{sec.cyclic-chern} 
we present the operator theoretic framework for the local index formula of~\cite{CM:LIFNCG}. Then we compute 
the CM cocycle of Dirac spectral triples: the even case is treated  in Section~\ref{sec.even-chern}  and the odd case in 
Section~\ref{sec.odd-chern}.

\section{Greiner's approach of the heat kernel asymptotics}
\label{sec.volterra}
  In this section we recall Greiner's approach of the heat kernel asymptotics as in \cite{Gr:AEHE} and 
\cite{BGS:HECRM} (see also~\cite[pp.~252-272]{Me:APSIT} for an alternative point of view). 

Here $M^n$  is a manifold equipped with a smooth and strictly positive density,    $\cE$  a Hermitian 
vector bundle over $M$ and  $\Delta$ a  second order elliptic differential operator on  $M$ acting on the sections 
of $\cE$. In addition we assume that $\Delta$ with domain 
$C^\infty_{c}(M,\cE)$ is essentially selfadjoint and bounded  from below on $L^{2}(M,\cE)$. Then by standard functional calculus 
we can define 
 $e^{-t\Delta}$, $t\geq0$, as a selfadjoint bounded 
operator on $L^{2}(M,\cE)$. 
In fact, $e^{-t\Delta}$ is smoothing for 
$t>0$ and so its distribution  kernel $k_{t}(x,y)$ belongs to 
$C^\infty(M,\cE)\hotimes C^\infty(M,\cE^{*}\otimes |\Lambda|(M))$ 
where $|\Lambda|(M)$ denotes the bundle of densities on $M$. 

Recall that the heat semigroup allows us to invert the heat equation, in the sense that the operator 
\begin{equation}
    Q_{0}u(x,s)=\int_{0}^\infty e^{-s\Delta} u(x,t-s)dt, \qquad u \in C^\infty_{c}(M\times \R, \cE), 
     \label{eq:volterra.inverse-heat-operator}
\end{equation}
maps continuously into $C^{0}(\R, L^{2}(M,\cE)) \subset \cD'(M\times\R, \cE)$ and satisfies
\begin{equation}
    (\Delta+\partial_{t})Q_{0}u = Q_{0}(\Delta+\partial_{t})u=u \qquad \forall u \in C^\infty_{c}(M\times\R,\cE).
\end{equation}

Notice that the operator $Q_{0}$  has the \emph{Volterra property} in the sense of~\cite{Pi:COPDTV}, i.e. it has a 
distribution kernel of the form $K_{Q_{0}}(x,y,t-s)$ where $K_{Q_{0}}(x,y,t)$ vanishes on the region $t<0$. In fact, 
\begin{equation}
    K_{Q_{0}}(x,y,t) = \left\{ 
    \begin{array}{ll}
         k_{t}(x,y) & \quad \text{if $t> 0$},  \\
        0 &  \quad \text{if $t<0$}. 
    \end{array}\right. 
\end{equation}

These equalities are the main motivation for using pseudodifferential techniques to study the heat kernel $k_{t}(x,y)$. 
The idea, which goes back to Hadamard~\cite{Ha:LCPLPDE}, is to consider a class of \psido's, the Volterra \psido's 
(\cite{Gr:AEHE},~\cite{Pi:COPDTV}, \cite{BGS:HECRM}), 
taking into account:  \smallskip 

(i)  The aforementioned Volterra property;   \smallskip

(ii) The parabolic homogeneity of the heat operator $\Delta+ \partial_{t}$, i.e. the homogeneity with respect 
             to the dilations $\lambda.(\xi,\tau)=(\lambda\xi,\lambda^{2}\tau)$, $(\xi,\tau)\in \R^{n+1}$, $\lambda\neq 0$. 
             \smallskip  
             
    In the sequel for $g\in \cS'(\R^{n+1})$  and $\lambda\neq 0$ we let $g_{\lambda}$ be the tempered distribution 
    defined by   
    \begin{equation}
        \acou{g_{\lambda}( \xi,\tau)}{u(\xi,\tau)} =    |\lambda|^{-(n+2)} 
            \acou{g(\xi,\tau)} {u(\lambda^{-1}\xi, \lambda^{-2}\tau)}, \quad u \in \cS(\R^{n+1}). 
    \end{equation}
\begin{definition}%
    A distribution $ g\in \cS'(\R^{n+1})$ is parabolic 
homogeneous of degree $m$, $m\in \Z$, if for any $\lambda \neq 0$ we have 
$g_{\lambda}=\lambda^m g$.  
\end{definition}
 
Let $\C_{-}$ denote the complex halfplane $\{\Im \tau >0\}$ with closure $\bar\C_{-}$. Then:  
\begin{lemma}[{\cite[Prop.~1.9]{BGS:HECRM}}] \label{lem:volterra.volterra-extension}
Let $q(\xi,\tau)\in C^\infty((\R^{n}\times\R)\setminus0)$ be a parabolic homogeneous 
symbol of degree $m$ such that: \smallskip 

(i) $q$ extends to a continuous function on $(\R^{n}\times\bar\C_{-})\setminus0$ in 
    such way to be holomorphic in the last variable when the latter is restricted to $\C_{-}$. \smallskip 

\noindent Then  there is a unique $g\in \cS'(\R^{n+1})$ agreeing with $q$ 
on $\R^{n+1}\setminus 0$ so that: \smallskip 

(ii) $g$ is homogeneous of degree $m$; \smallskip 

(iii) The inverse Fourier transform $\check g(x,t)$ vanishes for $t<0$. 
\end{lemma}

\begin{remark}%
    If we take $m\leq -(n+2)$ the result fails in general for symbols not satisfying~(i). 
\end{remark}
 
Let $U$ be an open subset of $\R^{n}$. We define Volterra symbols and Volterra \psido's on 
$U\times\R^{n+1}\setminus 0$ as follows. 

\begin{definition}%
    $S_{\op v}^m(U\times\R^{n+1})$, $m\in\Z$,  consists in smooth functions $q(x,\xi,\tau)$ on 
    $U\times\R^n\times\R$ with an asymptotic expansion  $q \sim \sum_{j\geq 0} q_{m-j}$ where: \smallskip 
    
    - $q_{l}\in C^{\infty}(U\times[(\R^n\times\R)\setminus0])$ is a homogeneous Volterra symbol of degree~$l$, 
    i.e. $q_{l}$ is  parabolic homogeneous of degree $l$ and  satisfies the 
    property (i) in Lemma~\ref{lem:volterra.volterra-extension} with respect to the last $n+1$ variables; \smallskip 
    
    - The sign $\sim$ means that, for any integer $N$ and any compact $K\subset U$, there is a constant $C_{NK\alpha\beta 
    k}>0$ such that 
            \begin{equation}
                |\partial^{\alpha}_{x}\partial^{\beta}_{\xi} \partial^k_{\tau}(q-\sum_{j< N} 
            q_{m-j})(x,\xi,\tau) | 
                \leq C_{NK\alpha\beta k} (|\xi|+|\tau|^{1/2})^{m-N-|\beta|-2k}, \qquad  
                          \label{eq:volterra.asymptotic-symbols}
            \end{equation}
    for  $x\in K$ and $|\xi|+|\tau|^{\frac12}>1$. 
\end{definition}

\begin{definition}\label{def:volterra.PsiDO}
    $\pvdo^m(U\times\R)$, $m\in\Z$,  consists in continuous operators 
    $Q$ from $C_{c}^\infty(U_{x}\times\R_{t})$ to 
    $C^\infty(U_{x}\times\R_{t})$ such that: \smallskip 
    
    (i) $Q$ has the Volterra property; \smallskip 
    
    (ii) $Q=q(x,D_{x},D_{t})+R$ for some symbol $q$ in $S^m_{\op v}(U\times\R)$ and some smoothing operator  $R$. 
\end{definition}

In the sequel if $Q$ is a Volterra \psido\ we let $K_{Q}(x,y,t-s)$ denote its distribution kernel, so that the 
distribution $K_{Q}(x,y,t)$ 
vanishes for $t<0$.  


\begin{example}%
    Let $P$ be a differential operator of order $2$ on $U$ and let $p_{2}(x,\xi)$ denote the principal symbol of $P$. Then 
    the heat operator $P+\partial_{t}$ is a Volterra \psido\ of order $2$ with principal symbol $p_{2}(x,\xi)+i\tau$. 
\end{example}

Other examples of Volterra \psido's are given by the homogeneous operators as in below.
\begin{definition}\label{def:volterra.homogeneous-PsiDO} 
Let $q_{m}(x,\xi,\tau) \in C^\infty(U\times(\R^{n+1}\setminus 0))$ be a homogeneous Volterra 
symbol of order $m$ and let $g_{m}\in C^\infty(U)\hotimes 
    \cS'(\R^{n+1})$ denote its unique homogeneous extension given by 
    Lemma~\ref{lem:volterra.volterra-extension}. Then: \smallskip 
    
    -  $\check q_{m}(x,y,t)$ is the inverse Fourier transform of $g_{m}(x,\xi,\tau)$ 
	in the last $n+1$ variables; \smallskip 
    
        -  $q_{m}(x,D_{x},D_{t})$ is the operator with kernel $\check q_{m}(x,y-x,t)$.
\end{definition}

\begin{proposition}[\cite{Gr:AEHE},~\cite{Pi:COPDTV}, \cite{BGS:HECRM}]\label{prop:Volterra-properties}
  The  following properties hold.\smallskip 

 \emph{1) Composition}. Let $Q_{j}\in \pvdo^{m_{j}}(\UR)$, $j=1,2$, have symbol $q_{j}$ 
and 
    suppose that $Q_{1}$ or $Q_{2}$ is properly supported. Then $Q_{1}Q_{2}$ belongs  to $\pvdo^{m_{1}+m_{2}}(\UR)$
    and has symbol  $q_{1}\#q_{2} \sim \sum \frac{1}{\alpha!} 
\partial_{\xi}^{\alpha}q_{1} D_{\xi}^\alpha q_{2}$.\smallskip 

   \emph{2) Parametrices.} An operator $Q\in \pvdo^{m}(U\times\R)$ admits a  parametrix 
    in $\pvdo^{-m}(U\times\R)$ if, and only if, its principal symbol is nowhere vanishing on 
    $U\times[(\R^n\times \bar\C_{-}\setminus 0)]$. \smallskip 

  \emph{3)  Invariance.} Let $\phi : U \rightarrow V$ be 
    a diffeomorphism onto another open 
subset $V$ of $\R^n$ and let $Q$ be a Volterra \pdo\ on $U\times\R$ of order $m$. Then 
$Q=(\phi\oplus \op{id}_{\R})_{*}Q$  is a Volterra \pdo\ on $V\times \R$ of order $m$. 
\end{proposition}

In addition to the above standard properties there is the one below which shows the relevance of 
Volterra \psido's for deriving small times asymptotics.  
\begin{lemma}[{\cite[Chap.~I]{Gr:AEHE}}, {\cite[Thm.~4.5] {BGS:HECRM}}]\label{lem:volterra.key-lemma}
Let $Q\in \pvdo^{m}(U\times\R)$ have symbol $q \sim \sum q_{m-j}$. Then the following asymptotics 
holds in  $C^\infty(U)$, 
\begin{equation}
    K_{Q}(x,y,t) \sim_{t\rightarrow 0^{+}} t^{-(\frac{n}2+[\frac{m}2]+1)} \sum_{l\geq 0} t^l 
             \check{q}_{2[\frac{m}2]-2l}(x,0,1), \label{eq:volterra.asymptotics-Q}
\end{equation}
     where the notation $\check{q}_{k}$ has the same meaning 
     as in Definition~\ref{def:volterra.homogeneous-PsiDO}. 
\end{lemma}
\begin{proof}
   As the Fourier transform relates the decay at infinity to 
    the behavior at the origin of the Fourier transform the 
    distribution $\check{q} -\sum_{j\leq J} \check{q}_{m-j}$ lies in $C^{N}
    (U_{x}\times\R^{n}_{y}\times\R_{t})$ 
    as soon as $J$ is large enough. Since $Q-q(x,D_{x},D_{t})$ is smoothing it  follows that 
   $R_{J}(x,t)=K_{Q}(x,x,t) -\sum_{j\leq J} \check q_{m-j}(x,0,t)$ is of class $C^{N}$. 
   As $R_{J}(x,y,t)=0$ for $t<0$ we get $\partial_{t}^l R_{J}(x,0)=0$ for 
    $l=0,1,\ldots,N$, so that $R_{J}(.,t)$ is a $\op{O}(t^{N})$ in $C^N(U)$ as $t\rightarrow 
    0^{+}$. It follows that in $C^\infty(U)$ we have the asymptotics  
        $K_{Q}(x,x,t) \sim_{t\rightarrow 
        0^{+}}   \sum \check{q}_{m-j}(x,0,t)$. 
        
    Now,  $(\check{q}_{m-j})_{\lambda}= |\lambda|^{-(n+2)} (q_{m-j,\lambda^{-1}})^{\vee}=  |\lambda|^{-(n+2)} \lambda^{j-m} 
    \check{q}_{m-j}$ for any $\lambda\neq 0$. So letting $\lambda=\sqrt{t}$, $t>0$, yields
 $\check  q_{m-j}(x,0,t)=t^{\frac{j-n-m}2-1} \check q_{m-j}(x,0,1)$, while for $\lambda=-1$ we get  
  $\check q_{m-j}(x,0,1)=-q_{m-j}(x,0,1)=0$ whenever $m-j$ is odd.  Thus,  
\begin{equation}
    K_{Q}(x,x,t)\sim_{t\rightarrow 0^{+}}  \sum_{m-j\ \text{even}} t^{\frac{j-n-m}2-1} \check{q}_{m-j}(x,0,1),
\end{equation}
 that is 
$K_{Q}(x,x,t) \sim_{t\rightarrow 0^{+}}  
    t^{-(\frac{n}2 + [\frac{m}2]+1)} \sum_{l\geq 0} t^l \check{q}_{2[\frac{m}2]-l}(x,0,1)$. 
  \end{proof}

The invariance property in Proposition~\ref{prop:Volterra-properties}
allows us to define Volterra \pdo's on $M\times \R$ acting on the 
sections of the vector bundle $\cE$. Then all the preceding  properties hold \emph{verbatim} in this context. 
In particular the heat operator $\Delta+\partial_{t}$  has a  parametrix  $Q$ in $\pvdo^{-2}(M,\times\R,\cE)$. 
In fact, comparing the operator~(\ref{eq:volterra.inverse-heat-operator}) with any 
Volterra parametrix for $\Delta+\partial_{t}$ allows us to prove:  
\begin{theorem}[\cite{Gr:AEHE},~\cite{Pi:COPDTV}, {\cite[pp.~363-362]{BGS:HECRM}}]\label{thm:volterra.inverse} 
    The differential operator $\Delta+\partial_{t}$ is invertible and its inverse 
    $(\Delta+\partial_{t})^{-1}$ is a Volterra \pdo\ of order $-2$. 
\end{theorem}

Combining this with Lemma~\ref{lem:volterra.key-lemma} gives the heat kernel asymptotics below.  
 
\begin{theorem}[{\cite[Thm.~1.6.1]{Gr:AEHE}}]\label{thm:volterra.heat-kernel-asymptotic}
  In  $C^\infty(M,|\Lambda|(M)\otimes \End\cE)$ we have:  
   \begin{equation}
        k_{t}(x,x) \sim_{t\rightarrow 0^{+}} t^{\frac{-n}2} \sum_{l\geq 0} t^l a_{l}(\Delta)(x), \qquad 
	a_{l}(\Delta)(x)=\check{q}_{-2-2l}(x,0,1), 
         \label{eq:volterra.heat-kernel-asymptotic}
    \end{equation}
where the 
equality on the right-hand side shows how to compute the densities $a_{l}(\Delta)(x)$'s in local trivializing coordinates 
by means of the symbol  $q\sim \sum q_{-2-j}$ of any Volterra parametrix for $\Delta+\partial_{t}$. 
\end{theorem}

This approach to the heat kernel asymptotics present several advantages. First, as 
Theorem~\ref{thm:volterra.heat-kernel-asymptotic} is a purely local statement we can easily localize the heat kernel 
asymptotics. In fact, given a Volterra parametrix $Q$ for $\Delta+\partial_{t}$ in some local trivializing coordinates around 
$x_{0}\in M$, comparing the asymptotics~(\ref{eq:volterra.asymptotics-Q}) and~(\ref{eq:volterra.heat-kernel-asymptotic}) we get
\begin{equation}
         k_{t}(x_{0},x_{0})=K_{Q}(x_{0},x_{0},t) +\op{O}(t^{\infty})\qquad  \text{as  $t\rightarrow 0^{+}$}.
         \label{eq:Volterra.kt-KQx0}
\end{equation}
    Therefore in order to determine the heat kernel asymptotics~(\ref{eq:volterra.heat-kernel-asymptotic}) at 
    $x_{0}$ we only need a  Volterra parametrix for 
    $\Delta+\partial_{t}$ near $x_{0}$.   

Second, we have a genuine asymptotics with respect to the $C^{\infty}$-topology, which can be differentiated as follows.
\begin{proposition}\label{prop:Volterra.asymptotics-Pheat}
    Let $P:C^{\infty}(M,\cE)\rightarrow C^{\infty}(M,\cE)$ be a differential operator of order $m$ and let $h_{t}(x,y)$ 
    denote the distribution kernel of $Pe^{-t\Delta}$. Then in $C^{\infty}(M,|\Lambda|\otimes \End \cE)$ we have 
    \begin{equation}
        h_{t}(x,x)\sim_{t\rightarrow 0^{+}} t^{[\frac{m}{2}]-\frac{n}{2}} \sum_{l \geq 0} t^{l} b_{l}(x), \qquad 
        b_{l}(x)= \check{r}_{2[\frac{m}{2}]-2-2l}(x,0,1),
         \label{eq:volterra.heat-kernel-asymptotic2}
    \end{equation}
    where the 
equality on the right-hand side gives a formula for computing the densities $b_{l}(x)$'s in local trivializing coordinates 
using the symbol  $r\sim \sum r_{m-2-j}$ of  $R=P(\Delta+\partial_{t})^{-1}$ (or of $R=PQ$ where $Q$ is any Volterra parametrix 
for $\Delta+\partial_{t}$). 
\end{proposition}
\begin{proof}
    As $h_{t}(x,y)=P_{x}k_{t}(x,y)=P_{x}K_{(\Delta+\partial_{t})^{-1}}(x,y,t)=K_{P(\Delta+\partial_{t})^{-1}}(x,y,t)$ 
    the result follows by applying Lemma~\ref{lem:volterra.key-lemma} to $P(\Delta+\partial_{t})^{-1}$ (or to
    $PQ$ where $Q$ is any Volterra parametrix for $\Delta+\partial_{t}$). 
\end{proof}

Finally, in local trivializing coordinates the densities $a_{j}(\Delta)(x)$'s can be explicitly computed in terms of the symbol  
    $p=p_{2}+p_{1}+p_{0}$ of $\Delta$.  To see this  
   let   $q\sim \sum q_{-2-j}$ be the symbol of a Volterra 
    parametrix $Q$ for $\Delta+\partial_{t}$. As  
    $q\#p\sim q(p+i\tau) +  \sum \frac{1}{\alpha!} 
     \partial_{\xi}^{\alpha}q D_{x}^\alpha p \sim 1$ we get $q_{-2}=(p_{2}+i\tau)^{-1}$ and 
    \begin{equation}
            q_{-2-j}= -(\!\sum_{k+l+|\alpha|=j} \!
            \frac1{\alpha!} \partial_{\xi}^\alpha q_{-2-k}   D_{x}^\alpha p_{2-l})(p_{2}+i\tau)^{-1}, \quad j\geq 1.
             \label{eq:volterra.symbols-parametrix}
       \end{equation}    
       Therefore, combining with~(\ref{eq:volterra.heat-kernel-asymptotic}) 
       we deduce that, as in~\cite{Gi:ITHEASIT}, the densities $a_{j}(\Delta)(x)$'s are 
       universal polynomials in the the jets at $x_{0}$ of the symbol of $\Delta$ with coefficients depending smoothly 
       on its principal symbol. Similarly, in local trivializing coordinates the densities $b_{l}(x)$'s 
       in~(\ref{eq:volterra.heat-kernel-asymptotic2}) can 
       be expressed  universal polynomials in the the jets at $x_{0}$ of the symbols of $\Delta$ and $P$ 
       with coefficients depending smoothly on the principal symbol of $\Delta$. 

\section{The local index formula of Atiyah and Singer}
\label{sec.getzler}
 In this section we shall give a new proof the local index formula of Atiyah and 
Singer~(\cite{AS:IEO1},~\cite{AS:IEO3}) by using Greiner's approach of the heat kernel asymptotics. 

Let $(M^n, g)$ be an even dimensional compact Riemannian spin manifold with spin bundle $\sS$
and let $\cE$ denote a Hermitian vector bundle over 
$M$ equipped with an unitary connection $\nabla^\cE$ 
with curvature $F^{\cE}$. Since $n$ is even  $\End \sS$ is as a bundle of 
algebras over $M$ isomorphic to the Clifford bundle $\Cl(M)$, whose  fiber $\Cl_{x}(M)$ at 
$x\in M$ is the complex algebra generated by $1$ and $T^{*}_{x}M$ with relations
\begin{equation}
    \xi.\eta + \eta.\xi = -2\acou\xi\eta,  \qquad \xi,\eta\in T^{*}_{x}M.
\end{equation}

Recall that he quantization map $c:\Lambda T^{*}_{\C}M \rightarrow \Cl(M)$ and the symbol map 
$\sigma=c^{-1}$ satisfy 
\begin{equation}
     \sigma(c(\xi)c(\eta)) = \xi\wedge\eta -\xi \llcorner  \eta, \qquad \xi \in T^{*}_{\C}M, \quad \eta \in \Lambda T^{*}_{\C}M  ,     
\end{equation}
where $\llcorner$ is the interior product.  Therefore, for  $\xi$ and $\eta$ in $\Lambda T^{*}_{\C}M$ we have   
\begin{equation}
    \sigma(c(\xi^{(i)})c(\eta^{(j)}))= \xi^{(i)}\wedge\eta^{(j)} \quad \bmod 
\Lambda^{i+j-2}T^{*}_{\C}M ,
    \label{eq:getzler.clifford-multiplication}
\end{equation}
where $\zeta^{(l)}$ denotes the component in $\Lambda^{l}T^{*}_{\C}M$ of $\zeta \in \Lambda T^{*}_{\C}M$. 
Thus the $\Z_{2}$-grading on $\Lambda T^{*}_{\C}M$ given by the parity of forms induces a $\Z_{2}$-grading 
 $\sS=\sS^{+}\otimes\sS^{-}$ on the spin bundle. Furthermore, if $e_{1},\ldots,e_{n}$ is an 
 orthonormal frame for $T_{x}M$ and  we regard $c(dx^{i_{1}})\cdots 
c(dx^{i_{k}})$, $i_{1}<\ldots <i_{k}$, as an endomorphism of $\sS_{x}$ then 
\begin{equation}
   \Str_{x} c(e^{i_{1}})\cdots c(e^{i_{k}}) =\left\{  \begin{array}{cl}
        0 & \text{if $k\neq n$,}   \\
  (-2i)^{\frac{n}2}  &\text{if $k=n$}.
    \end{array} \right.
    \label{eq:asymptotic.supertrace-even}
\end{equation}

Let  $\nabla^{\sSE}=\nabla^\sS\otimes 1 + 1\otimes 
    \nabla^\cE$ be the connection on $\sSE$, where $\nabla^\sS$ denotes the Levi-Civita connection lifted to the spin bundle. 
    Then the Dirac operator $\sD_{\cE}$ acting on the sections of $\sSE$ 
is given by the composition 
\begin{equation}
    C^\infty(M, \sSE) \stackrel{\nabla^{\sSE}}{\longrightarrow} C^\infty(M,T^{*}M \otimes\sSE) 
    \stackrel{c\otimes 1}{\longrightarrow} C^\infty(M,\sSE).
\end{equation}
It is odd with respect to the $\Z_{2}$-grading 
    $\sSE=(\sS^+\otimes\cE)\oplus(\sS^-\otimes\cE)$, i.e. it can be written in the form 
\begin{equation}
            \sD_{\cE} = \left( 
            \begin{array}{cc}
                0 & D_{\cE}^{+}  \\ 
                D_{\cE}^{-}& 0
            \end{array}
            \right) , \qquad \sD_{\cE}^{\pm}: C^{\infty}(M, \sS^\mp\otimes\cE) \rightarrow C^{\infty}(M, \sS^\pm\otimes\cE). 
\end{equation}
Moreover, by the Lichnerowicz formula (\cite{BGV:HKDO},~\cite{LM:SG}, \cite{Roe:EOTAM}) we have 
 \begin{equation}
 \sD^{2}_{\cE}=
(\nabla^{\sSE}_{i})^{*}\nabla^{\sSE}_{i} + \mathcal{F^\cE}
 + \frac{\kappa^{M}}4,    
    \label{eq:AS.lichnerowicz}
\end{equation}
where $\kappa^{M}$ denotes the scalar curvature of $M$ and $\mathcal{F^\cE}$  the curvature $F^\cE$ 
lifted to $\sSE$, i.e. $F^\cE=\frac{1}{2}c(e^{k})c(e^{l})F^{\cE}(e_{k},e_{l})$ for any local orthonormal tangent frame 
$e_{1},\ldots,e_{n}$.   It follows that  $\sD_{\cE}$ and $\sD_{\cE}^{\pm}$ are elliptic, hence   
 are Fredholm.

 \begin{theorem}[\cite{AS:IEO1},~\cite{AS:IEO3}]\label{thm:aymptotic-AS} We have: 
   \begin{equation}
   \ind \sD_{\cE}^{+} = (2i\pi)^{-\frac{n}2} \int_{M} [\hat A(R_{M})\wedge \Ch (F^\cE)]^{(n)}, 
	\label{eq:asymptotic-AS-formula}
    \end{equation}
where $\hat A(R^{M})=\det{}^{\frac12}( \frac{R^{M}/2}{\sinh (R^{M}/2)}) $ is the total $\hat A$-form of the Riemann curvature and 
    $\Ch(F^\cE)=\Tr  \exp(-F^\cE)$ the total Chern form of the curvature $F^{\cE}$. 
\end{theorem}

In fact,  by the McKean-Singer formula $\ind \sD_{\cE}^{+}  =  \Str e^{-t\sD_{\cE}^{2}}$ for any  $t>0$.  
Therefore the index formula  follows from: 
\begin{theorem}\label{thm:index.heat-kernel-asymptotic}
In $C^\infty(M, |\Lambda|(M))$ we have  
\begin{equation}
    \Str_{x} k_{t}(x,x) = [\hat A(R_{M})\wedge \Ch (F^\cE)]^{(n)} + \op{O}(t) \qquad \text{as $t\rightarrow 0^{+}$}. 
     \label{eq:AS.local-index-theorem}
\end{equation}
\end{theorem}

This theorem, also called local index theorem,  was 
first proved by Patodi, Gilkey and 
Atiyah-Bott-Patodi~(\cite{ABP:OHEIT}, \cite{Gi:ITHEASIT}), and then in a purely analytic fashion by 
Getzler~(\cite{Ge:POSASIT},~\cite{Ge:SPLASIT}) and Bismut~\cite{Bi:ASITPA} 
(see also~\cite{BGV:HKDO}, \cite{Roe:EOTAM}). 
Moreover, as it is a purely local statement it holds \emph{verbatim} for 
(geometric) Dirac operators acting on a Clifford bundle. Thus it allows us to recover, on the one hand, the 
Gauss-Bonnet, signature and Riemann-Roch theorems (\cite{ABP:OHEIT}, \cite{BGV:HKDO}, \cite{LM:SG}, 
\cite{Roe:EOTAM})  and, on the other hand, the full 
index theorem of Atiyah-Singer (\cite{ABP:OHEIT}, \cite{LM:SG}). 

The short proof of  Getzler~\cite{Ge:SPLASIT} combines the Feynman-Kac representation of the heat kernel with an ingenious trick, 
the Getzler rescaling. We can alternatively prove Theorem~\ref{thm:index.heat-kernel-asymptotic} by combining Getzler 
rescaling with Greiner's approach of the heat kernel asymptotics as follows. 
\begin{proof}[Proof of Theorem~\ref{thm:index.heat-kernel-asymptotic}]
First, the Greiner approach allows us to easily localize the problem (compare~\cite{Ge:SPLASIT}).  Indeed, thanks to  
Theorem~\ref{thm:volterra.heat-kernel-asymptotic}  $ \Str_{x} k_{t}(x,x)$ admits an asymptotics 
in $C^\infty(M, |\Lambda|(M))$ as $t\rightarrow 0^{+}$. Thus, it is enough to prove~(\ref{eq:AS.local-index-theorem}) 
at a point $x_{0}\in M$. 
Furthermore, to reach this aim we know from~(\ref{eq:Volterra.kt-KQx0}) that we only need a Volterra parametrix for 
$\sD_{\cE}^{2}+\partial_{t}$ in local trivializing coordinates centered at $x_{0}$. Therefore, using normal coordinates 
centered at $x_{0}$ and a trivialization of the tangent bundle by means of a synchronous frame 
$e_{1},\ldots, e_{n}$ such that $e_{j}=\partial_{j}$ at $x=0$ we may replace $\sD_{\cE}$ by a Dirac operator $\sD$ on $\R^n$ 
acting on the trivial bundle with fiber 
$\sS_{n}\otimes\C^p$, where $\sS_{n}$ denotes the spin bundle of $\R^n$. Then we have
\begin{equation}
    k_{t}(0,0)=K_{Q}(0,0,t)+\op{O}(t^{\infty}) \qquad \text{as $t\rightarrow 0^{+}$}.
    \label{eq:AS.kt-KQ0}
\end{equation}

Second, as pointed out in~\cite{ABP:OHEIT} (see 
also~\cite{BGV:HKDO}, \cite{Roe:EOTAM}) choosing normal coordinates and a synchronous tangent frame makes  
the metric $g$ and the coefficients $\omega_{ikl}= \acou{\nabla^{LC}_{i}e_{k}}{e_{l}}$ of the 
Levi-Civita connection have behaviors near $x=0$ of the form  
\begin{equation}
     g_{ij}(x)=\delta_{ij}+\op{O}(|x|^{2}), \qquad  \omega_{ikl}(x)= -\frac12 
R_{ijkl}^{M}(0)x^j 
     +\op{O}(|x|^{2}),   
    \label{eq:AS-asymptotic-geometric-data}
\end{equation}
where  $R_{ijkl}^{M}(0)=\acou{R^{M}(0)(\partial_{i},\partial_{j})\partial_{k}}{\partial_{l}}$.  
Then using~(\ref{eq:asymptotic.supertrace-even}) and~(\ref{eq:AS.kt-KQ0}) we get
\begin{equation}
    \Str k_{t}(0,0) =(-2i)^{\frac{n}2} \sigma\otimes \Tr_{\C^{p}}[K_{Q}(0,0,t)]^{(n)}+ 
\op{O}(t^{\infty}). 
     \label{eq:AS.kt-sKQ}
\end{equation}
Thus we are reduced  to prove the convergence of $\sigma[K_{Q}(0,0,t)]^{(n)}$ as $t\rightarrow 0^{+}$ and to identify its limit. 

Now, recall that the Getzler rescaling~\cite{Ge:SPLASIT} assigns the following degrees:
\begin{equation}
  \deg \partial_{j}=\frac12 \deg \partial_{t}=\deg c(dx^j)=- \deg x^j =1 , \qquad 
  \label{eq:AS-Getzler-order}
\end{equation}
while $ \deg B=0$ for any $B\in M_{p}(\C)$. 
It can define a filtration of Volterra \psido's with coefficients in $\End (\sS_{n}\otimes \C^{p}) 
\simeq \Cl(\R^{n})\otimes M_{p}(\C)$ as follows.  

Let $Q\in  \Psi_{\op{v}}^{*}(\R^{n}\times\R, \sS_{n}\otimes \C^{p})$ have symbol $q(x,\xi,\tau) \sim \sum_{k\leq m'} 
q_{k}(x,\xi,\tau)$. Then taking components in each subspace 
 $\Lambda^j T_{\C}^{*}\R^{n}(n)$ and then using Taylor expansions at $x=0$ gives formal expansions  
\begin{equation}
    \sigma[q(x,\xi,\tau)] \sim \sum_{j,k} \sigma[q_{k}(x,\xi,\tau)]^{(j)} \sim \sum_{j,k,\alpha} 
    \frac{x^{\alpha}}{\alpha!}  \sigma[\partial_{x}^{\alpha}q_{k}(0,\xi,\tau)]^{(j)}.
    \label{eq:AS.Getzler-asymptotic}
\end{equation}
According to~(\ref{eq:AS-Getzler-order}) the symbol $ \frac{x^{\alpha}}{\alpha!}  
\partial_{x}^{\alpha}\sigma[q_{k}(0,\xi,\tau)]^{(j)}$ is Getzler homogeneous of degree $k+j-|\alpha|$. Therefore,  we can expand 
$\sigma[q(x,\xi,\tau)]$ as 
\begin{equation}
    \sigma[q(x,\xi,\tau)] \sim \sum_{j\geq 0} q_{(m-j)}(x,\xi,\tau), \qquad q_{(m)}\neq 0, 
    \label{eq:index.asymptotic-symbol}
\end{equation}
where $q_{(m-j)}$ is a Getzler homogeneous symbol of degree $m-j$.

\begin{definition}\label{def:getzler.model-operator}
Using~(\ref{eq:index.asymptotic-symbol}) we make the following definitions: 
 \smallskip 
  
  - The integer $m$   is the Getzler order of $Q$,\smallskip 
  
  - The symbol $q_{(m)}$ is the principal Getzler homogeneous symbol  of $Q$, \smallskip 
  
  - The operator $Q_{(m)}=q_{(m)}(x,D_{x},D_{t})$ is the model operator of $Q$.
\end{definition}

\begin{remark}%
 The model operator $Q_{(m)}$ is well defined according to definition~\ref{def:volterra.homogeneous-PsiDO}. \smallskip 
\end{remark}

\begin{remark}%
 By construction we always have $\text{Getzler order}\leq \text{order} + n$, but this is not  an equality in general. 
\end{remark}

\begin{example}%
Let $A=A_{i}dx^i$ is the connection one-form on $\C^p$. Then  by~(\ref{eq:AS-asymptotic-geometric-data}) the covariant 
derivative $ \nabla_{i}= \partial_{i}+ \frac14 \omega_{ikl}(x)c(e^k)c(e^l) +A_{i}$  on $\sS_{n}\otimes \C^p$
    has Getzler order 1 and model operator  
    \begin{equation}
        \nabla_{i (1)} =\partial_{i}-\frac14 R_{ij}^{M}(0)x^j , \qquad 
        R_{ij}^{M}(0)=R_{ijkl}^{M}(0)dx^k \wedge dx^l. 
        \label{eq:AS.model-spin-connection}
    \end{equation}
\end{example}

The interest to introduce Getzler orders stems from the following. 
\begin{lemma}\label{lem:AS.approximation-asymptotic-kernel}
     Let $Q\in \Psi_{\op{v}}^{*}(\R^{n}\times\R, S_{n}\otimes \C^{p})$ have Getzler order $m$ and model operator 
     $Q_{(m)}$. Then 
      as  $t\rightarrow 0^{+}$ we have: \smallskip 
 
  -  $\sigma[K_{Q}(0,0,t)]^{(j)}= \op{O}(t^{\frac{j-m-n-1}2})$ if $m-j$ is odd; \smallskip 

         - $\sigma[K_{Q}(0,0,t)]^{(j)}= t^{\frac{j-m-n}2-1} 
         K_{Q_{(m)}}(0,0,1)^{(j)} + \op{O}(t^{\frac{j-m-n}2}) $ if  $m-j$ is even.   \smallskip 
  
\noindent In particular for $m=-2$ we get 
\begin{equation}
    \sigma[K_{Q}(0,0,t)]^{(n)}= K_{Q_{(-2)}}(0,0,1)^{(n)} + \op{O}(t). 
    \label{eq:AS-convergence-symbol-KQ}
\end{equation}
 \end{lemma}
\begin{proof}
 Let $q(x,\xi,\tau)\sim \sum q_{k}(x,\xi,\tau)$ be the symbol of $Q$ and let $q_{(m)}(x,\xi,\tau)$ be the principal 
 Getzler homogeneous symbol. By  
 Lemma~\ref{lem:volterra.key-lemma} we have 
\begin{equation}
     \sigma[K_{Q}(0,0,t)]^{(j)}\sim_{t\rightarrow 0^{+}}  \sum t^{-\frac{n+2+m-j}{2}} \sigma[\check{q}_{k}(0,0,1)]^{(j)},
\end{equation}
and we know that $\check{q}_{k}(0,0,1)=0$ if $k$ is odd. Also, the symbol $\sigma[q_{k}(0,\xi,\tau)]^{(j)}$ is 
Getzler homogeneous of degree $k+j$, so it must be zero if $k+j>m$ since otherwise  $Q$ would not have Getzler order $m$. 
Hence: 
\smallskip 
 
 -  $\sigma[K_{Q}(0,0,t)]^{(j)}= \op{O}(t^{\frac{j-m-n+1}2})$ if $m-j$ is odd; \smallskip  
         
  \indent - $\sigma[K_{Q}(0,0,t)]^{(j)}= t^{\frac{j-m-n}2-1} 
         \sigma[\check{q}_{m-j}(0,0,1)]^{(j)} + \op{O}(t^{\frac{j-m-n}2}) $ if  $m-j$ is even.\smallskip 

\noindent On the other hand, notice that the symbol 
 $ \sigma[q_{(m)}(0,\xi,\tau)]^{(j)}$ is equal to 
 \begin{equation}
     \sum_{k+j-|\alpha|=m}  (\frac{x^{\alpha}}{\alpha!}  
     \partial_{x}^{\alpha}\sigma[q_{k}(0,\xi,\tau)]^{(j)})_{x=0}= \sigma[q_{m-j}(0,\xi,\tau)]^{(j)}.
 \end{equation}
 Thus $\sigma[\check{q}_{m-j}(0,0,1)]^{(j)}=K_{Q_{(m)}}(0,0,1)^{(j)}$. Hence the lemma.
\end{proof}

In the sequel we say that a symbol or a \psido\ is $\op{O_{G}}(m)$ if it has Getzler order~$\leq m$. 
\begin{lemma}\label{lem:index.top-total-order-symbol-composition}
For  $j=1,2$ let $Q_{j}\in \Psi^{*}_{\op{v}}(\R^{n}\times\R, \End(\sS_{n}\otimes \C^{p})$ have Getzler order $m_{j}$ 
and model 
operator $Q_{(m_{j})}$ and assume either $Q_{1}$ or $Q_{2}$ properly supported. Then we have: 
\begin{equation}
    Q_{1}Q_{2}= c[Q_{(m_{1})} Q_{(m_{2})}] +\op{O}_{G}(m_{1}+m_{2}-1).
\end{equation}
\end{lemma}
\begin{proof}
Let $q_{j}$ be the symbol of $Q_{j}$ and let $q_{(m_{j})}$ be 
its principal Getzler homogeneous symbol. By Proposition~\ref{prop:Volterra-properties} the operator 
$Q_{1}Q_{2}$ has symbol $q_{1}\#q_{2}$. 
Moreover for $N$ large enough  $q_{1}\# q_{2} -  \sum_{|\alpha|\leq N}\frac1{\alpha!} 
\partial_{\xi}^\alpha q_{1}  . 
D_{x}^\alpha q_{2}$ has order $<m_{1}+m_{2}-n$, so has Getzler order $<m_{1}+m_{2}$. As 
$\partial_{\xi}^\alpha q_{1}.  D_{x}^\alpha q_{2} - c[\partial_{\xi}^\alpha q_{(m_{1})}\wedge 
D_{x}^\alpha  f_{(m_{2})}]$ has Getzler order $\leq m_{1}+m_{2}-|\alpha|-1$ it follows that for $N$ large enough,
\begin{equation}
    q_{1}\# q_{2} = \sum_{|\alpha|\leq N}\frac1{\alpha!} c(\partial_{\xi}^\alpha 
q_{m_{1}}  \wedge D_{x}^\alpha q_{m_{2}}) + \op{O}_{G}(m_{1}+ m_{2}-1).  
    \label{eq:AS.composition.symbol.eq2}
\end{equation} 

On the other hand, $\sum \frac1{\alpha!} \partial_{\xi}^\alpha 
q_{(m_{1})}  \wedge D_{x}^\alpha q_{(m_{2})}$ is exactly the symbol of $Q_{(m_{1})}Q_{(m_{2})}$ 
  since $q_{(m_{2})}(x,\xi,\tau)$ is polynomial in 
$x$ and thus the sum is finite. Therefore taking $N$ large enough 
in~(\ref{eq:AS.composition.symbol.eq2}) shows that the symbols of $Q_{1}Q_{2}$ and $Q_{(m_{1})}Q_{(m_{2})}$ 
coincide modulo a symbol of Getzler order $\leq m_{1}+m_{2}-1$.
\end{proof}

Recall that by the Lichnerowicz formula~(\ref{eq:AS.lichnerowicz}) we have 
\begin{equation}
    \sD^{2}_{\cE}= -g^{ij}(\nabla_{i}\nabla_{j} -\Gamma_{ij}^k 
     \nabla_{k}) 
        + \frac12 c(e^i)c(e^j)F(e_{i},e_{j})  + \frac{\kappa}4,
     \label{eq:AS.lichnerowicz-bis}
\end{equation}
where the $\Gamma_{ij}^k$'s are the Christoffel symbols of the metric. 
Thus combining Lemma~\ref{lem:index.top-total-order-symbol-composition} 
with~(\ref{eq:AS-asymptotic-geometric-data}) and~(\ref{eq:AS.model-spin-connection}) shows 
that $\sD^{2}$ has Getzler order $2$ and  its model operator is 
\begin{equation}
\begin{split}
\sD^{2}_{(2)}&= -\delta_{ij} \nabla_{i(1)}\nabla_{j(1)} + \frac12 
    F^\cE(\partial_{k},\partial_{l})(0)dx^k \wedge dx^l \\ 
    & = H_{R} +F^\cE(0),  \qquad H_{R}=- \sum_{i=1}^n (\partial_{i}-\frac14 
R_{ij}^{M}(0)x^j)^{2}. 
\end{split}
\end{equation}

\begin{lemma}\label{lem:getzler-model-asymptotic}
Let $Q$ be a Volterra parametrix for $\sD^{2}+\partial_{t}$. Then: \smallskip 

1) $Q$ has Getzler order $2$ and its model operator is 
$(H_{R}+F^\cE(0)+\partial_{t})^{-1}$. \smallskip 

2) We have 
\begin{equation}
    K_{(H_{R}+F^\cE(0)+\partial_{t})^{-1}}(x,0,t)=G_{R}(x,t)\wedge e^{-tF^{\cE}(0)},
\end{equation}
 where $G_{R}(x,t)$ is the fundamental solution of  $H_{R}+\partial_{t}$, i.e. 
     the unique distribution such that 
     $(H_{R}+F^\cE(0)+\partial_{t})G_{R}(x,t)=\delta(x,t)$. \smallskip 
     
    3) As $t \rightarrow 0^{+}$ we have 
    \begin{equation}
    \sigma [K_{Q}(0,0,t)]^{(2j)}   =  t^{j-\frac{n}2} [G_{R}(0,1)\wedge e^{-F^{\cE}(0)}]^{(2j)} + \op{O}(t^{j-\frac{n}2+1}).
	     \label{eq:AS.approximation-asymptotic-kernel}
     \end{equation}
\end{lemma}
\begin{proof}
Note that 3) follows by combining 1) and 2) with 
Lemma~\ref{lem:AS.approximation-asymptotic-kernel}, so we only have to prove the first two assertions.  
Let $p(x,\xi) =\sum p_{j}(x,\xi)$ be the symbol of $\sD^{2}$ and let 
$q\sim \sum q_{-2-j}$ denote the symbol 
of $Q$. As $\sD^{2}$ is elliptic and has Getzler order $2$ we have 
$p_{(2)}(0,\xi)^{(0)}=p_{2}(0,\xi)\neq 0$.  Hence $q_{-2}=(p_{2}+i\tau)^{-1}$ has Getzler order $-2$. 
It then follows from~(\ref{eq:volterra.symbols-parametrix})
that each symbol $q_{-2-j}$ has Getzler order $\leq -2$. Hence $Q$ 
has Getzler order $-2$. 

On the other hand, $(\sD^{2}+\partial_{t})Q-1$ is smoothing, so  
by Lemma~\ref{lem:index.top-total-order-symbol-composition} the operator
$(H_{R}+F^\cE(0)+\partial_{t})Q_{(-2)}-1$ has Getzler order $\leq -1$. 
As the latter is Getzler homogeneous of degree~$0$ it must be zero. Hence 
$Q_{(-2)}=(H_{R}+F^\cE(0)+\partial_{t})^{-1}$, so that we have 
\begin{equation}
(H_{R,x}+F^\cE(0)+\partial_{t})K_{Q_{(-2)}}(x,y,t-s)=\delta(x-y,t-s).   
\label{eq:AS.inverse-model-kernel}
\end{equation}
Now, setting $y=0$ and $s=0$ in~(\ref{eq:AS.inverse-model-kernel}) shows that $G_{R,F}(x,t)= K_{Q_{(-2)}}(x,0,t)$ is the 
fundamental solution of $H_{R}+F(0)+\partial_{t}$. 
In fact, if we let  $G_{R}(x,t)$ be the fundamental solution of $H_{R}+F(0)+\partial_{t}$ then 
$G_{R,F}(x,t)=G_{R}(x,t)\wedge e^{-tF^{\cE}(0)}$. Thus $K_{Q_{(-2)}}(x,0,t)=G_{R}(x,t)\wedge 
e^{-tF^{\cE}(0)}$. 
\end{proof}

At this stage remark that $H_{R}$ is the harmonic oscillator associated to the antisymmetric matrix 
$R^{M}(0)=(R^{M}_{ij}(0))$.  Therefore we can make use of a version of the Melher formula~(\cite{GJ:QPFIPV}, 
\cite{Ge:SPLASIT}) to obtain:
\begin{lemma}\label{lem.getzler.melher-formula}
The fundamental solution $G_{R}(x,t)$ of  $H_{R}+\partial_{t}$ is  
\begin{equation}
     \chi(t) 
   (4\pi t)^{-\frac{n}2}  \det{}^{\frac12}( \frac{tR^{M}(0)/2}{\sinh (tR^{M}(0)/2)})
   \exp(-\frac1{4t} \acou{\frac{tR^{M}(0)/2}{\tanh (tR^{M}(0)/2)}x}x),
\end{equation}
where $\chi(t)$ is the characteristic function of $(0,+\infty)$. 
 \end{lemma}
\begin{proof}
  Let $a \in \R$ and let $H_{a}$ denote the harmonic oscillator $-\frac{d}{dx^{2}} + 
    \frac{1}4 a^{2}x^{2}$ on $\R$. Then the fundamental solution of $H_{a}+\partial_{t}$ is 
    $G_{a}(x,t)=\chi(t)S_{a}(x,t)$, where
\begin{equation}
       S_{a}(x,t) = (4\pi t)^{-\frac12} (\frac{at}{\sinh at})^{\frac12} 
       \exp(-\frac{1}{4t} x^{2}\frac{at}{\tanh at}), \quad t>0.   
\end{equation}
In fact $(H+\partial_{t})S_{a}=0$ on 
       $\R\times(0,+\infty)$ and $S(.,t)\rightarrow \delta$ in $\cS'(\R)$, since on compact 
       sets $ \hat S_{x \rightarrow \xi}(\xi,t) =\cosh^{-\frac12}(at) \exp(-\xi^{2}t \, 
\frac{\tanh at}{at}) $ converges to $1$. Hence 
$(H+\partial_{t})k_{a}=\chi'G(.,0)+\chi (H+\partial_{t})G=\delta$. 
    
More generally, if $A$ is a real $n\times n$ antisymmetric matrix and we let $B=-A^{2}$, then the fundamental solution of 
$-\sum \partial_{j}^{2} +\frac14 B_{jk}x^jx^k 
+\partial_{t}$ on $\R^{n}\times \R$ is 
\begin{equation}
    G_{A}(x,t)= \chi(t) 
   (4\pi t)^{-\frac{n}2} \det{}^{\frac12}( \frac{iAt}{\sinh (iAt)})
   \exp(-\frac1{4t} \acou{\frac{iAt}{\tanh (iAt)}x}x). 
     \label{eq:getzler.Melher-matrix}
\end{equation}
The passage from the formula for $G_{a}$ to the one for $G_{A}$ uses $O(n)$-invariance and in 
particular invariance under rotations in the $(x_{j},x_{k})$-plane,  $j<k$. Thus $G_{A}$ is 
also the fundamental solution for $-\sum (\partial_{j}-\frac{i}2 A_{jk}x^j)^{2} 
+\partial_{t}$. 

Now, the r.h.s. in~(\ref{eq:getzler.Melher-matrix}) is analytic with respect to $A$ 
and $R^{M}(0)$ is an antisymmetric matrix  made out of 2-forms which  
commute with other forms. Therefore the formula for $G_{A}$ with $A$ replaced by $-iR^{M}(0)/2$ 
gives the fundamental solution of $H_{R}+\partial_{t}$.  
\end{proof}

Finally, combining the formula for $G_{R}(x,t)$ with 
Lemma~\ref{lem:getzler-model-asymptotic} and (\ref{eq:AS.kt-sKQ}) we get 
\begin{equation}
    \Str k_{t}(0,0) = (2i\pi)^{-\frac{n}2}[ \hat{A}(R^{M}(0))\wedge \Ch(F^\cE(0))]^{(n)} + \op{O}(t) \quad \text{as $t\rightarrow 
    0^{+}$}. 
\end{equation}
This completes the proof of Theorem~\ref{thm:index.heat-kernel-asymptotic} and of the Atiyah-Singer index formula. 
\end{proof}

The main new feature in the previous proof is the use of Lemma~\ref{lem:AS.approximation-asymptotic-kernel} 
which, by very elementary considerations on 
Getzler orders, shows that the convergence of the supertrace of the heat kernel is a consequence of a general fact about 
Volterra \psidos. It also gives a differentiable version of Theorem~\ref{thm:index.heat-kernel-asymptotic} as follows.

In the sequel we abbreviate by \emph{synchronous normal coordinates centered at} 
$x_{0}\in M$ the data of normal coordinates 
centered at $x_{0}$ and of a trivialization of the tangent bundle $TM$ by means of a synchronous frame as in the proof 
of Theorem~\ref{thm:index.heat-kernel-asymptotic}. 

\begin{definition}%
    We say that $Q\in \Psi_{\op{v}}^{*}(M\times \R, \sS\otimes \cE)$ has Getzler order $m$ if it has 
    Getzler order $m$ in synchronous normal coordinates centered at any $x_{0}\in M$. 
\end{definition}

\begin{proposition}\label{prop:AS.differantiated-local-index-theorem}
    Let $\sP$ be a differential operator on $M$ acting on $\sS\otimes \cE$ whose Getzler order is equal 
    to $m$ and let $h_{t}(x,y)$ denote the kernel of $\sP e^{-t\sD_{\cE}^{2}}$. Then as $t\rightarrow 0^{+}$ we 
    have an asymptotics in $C^{\infty}(M,  |\Lambda|(M))$ of the form: \smallskip 
 
        - $\Str_{x}h_{t}(x,x)= \op{O}(t^{\frac{-m+1}2})$ if $m$ is odd;  \smallskip  
  
        -  $\Str_{x}h_{t}(x,x)= t^{\frac{-m}2} 
        B_{0}(\sD^{2}_{\cE},\sP)(x) + \op{O}(t^{\frac{-m}2+1}) $ if  $m$ is even, where  in synchronous normal coordinates 
        centered at $x_{0}$ and with $\sP_{(m)}$ denoting the model operator of $\sP$  we have 
        $B_{0}(\sD^{2}_{\cE},\sP)(0)= (-2i\pi)^{\frac{n}2} [(\sP_{(m)}G_{R})(0,1)\wedge \Ch(F^{\cE}(0))]^{(n)}$. 
\end{proposition}
\begin{proof}
As in the proof of Proposition~\ref{prop:Volterra.asymptotics-Pheat} 
we have $ h_{t}(x,y)= K_{\sP(\sD^{2}_{\cE}+\partial_{t})^{-1}}(x,y,t)$. Notice 
that  by Lemma~\ref{lem:index.top-total-order-symbol-composition} and Lemma~\ref{lem:getzler-model-asymptotic} 
in  synchronous normal coordinates $\sP(\sD^{2}_{\cE}+\partial_{t})^{-1}$ 
has Getzler order $m-2$ and its model operator is $Q_{(m-2)}=\sP_{(m)}(H_{R}+F^{\cE}+\partial_{t})^{-1}$. Thus 
$K_{Q_{(m-2)}}(x,0,t)=\sP_{(m)x} K_{(H_{R}+F^{\cE}+\partial_{t})^{-1}}(x,0,t)
     =(\sP_{(m)}G_{R})(x,t)\wedge e^{-tF^{\cE}(0)}$. Then the proposition follows by applying 
     Proposition~\ref{prop:Volterra.asymptotics-Pheat}  and 
     Lemma~\ref{lem:AS.approximation-asymptotic-kernel}.
%
%
%
%
\end{proof}


\section{The local index formula in noncommutative geometry}
\label{sec.cyclic-chern}

 In this section we recall the operator theoretic framework for the local index formula (\cite{Co:NCG}, 
\cite{CM:LIFNCG}; see also~\cite{Hi:LIFNCG}). This uses two main tools, spectral triples and cyclic cohomology. 

A \emph{spectral triple} is a triple $(\cA, \cH, D)$ where the involutive unital algebra 
$\cA$ is represented in the (separable) Hilbert space $\cH$ and $D$ is an unbounded 
selfadjoint 
operator on $\cH$ with compact resolvent and which almost commutes with $\cA$, i.e. $[D,a]$ is 
bounded for any element $a$ of $\cA$. 

In the sequel we assume $\cA$ 
stable by holomorphic calculus, i.e. if $a \in \overline{\cA}$ 
is invertible then $a^{-1} \in \cA$; this has the effect that the $K$-groups of $\cA$ and $\bar \cA$ coincide.  

The spectral triple is \emph{even} if $\cH$ is endowed with  a $\Z_{2}$-grading $\gamma \in 
    \cL(\cH)$, $\gamma=\gamma^{*}$, $\gamma^{2}=1$, such that $ \gamma D=-D\gamma$ and 
    $\gamma a =a \gamma$ for all $a\in \cA$.  
Otherwise the spectral triple is \emph{odd}.    

The datum of $D$ above defines an additive index map $\ind_{D}:K_{*} \rightarrow \Z$ 
as follows (see also~\cite[sect. 2]{Mo:IEDPNG}). 

In the even case, with respect to the decomposition $\cH=\cH^{+}\oplus \cH^-$ given by the 
$\Z_{2}$-grading of $\cH$ the operator $D$ takes the form 
\begin{equation}
     D= \left( 
        \begin{array}{cc}
            0& D^-  \\ 
            D^+& 0
        \end{array}
        \right)    \qquad D_{\pm} :\cH^{\mp} \rightarrow \cH^{\pm}.  
\end{equation}
For any selfadjoint idempotent $e \in M_{q}(\cA)$ the operator 
    $e(D^{+}\otimes 1)e$ from $e(\cH^{+}\otimes \C^q)$ to $e(\cH^{-}\otimes \C^q)$ is 
    Fredholm and its index only depends on the homotopy class of $e$. We then define 
    \begin{equation}
        \ind_{D}[e]= \ind eD^{+}e. 
        \label{eq:LIFNCG.index-even}
    \end{equation}

   In the odd case,  given an unitary $U\in GL_{q}(\cA)$ the operator 
$[D,U]$ is bounded and so the compression $P U P$, where $P=\frac{1+F}2$ and $F=\op{sign}D$, is Fredholm. The  index of 
$P U P$ then depends only on the homotopy class of $U$ and we let 
\begin{equation}
    \ind_{D}[U]= \ind P U P. 
        \label{eq:LIFNCG.index-odd}
\end{equation}

The index map~(\ref{eq:LIFNCG.index-odd}) can also be interpreted in terms of spectral flows as follows. 
Recall that given a family $(D_{t})_{0\leq t\leq 1}$ of (unbounded) selfadjoint operators with 
discrete spectrum such that $D_{0}-D_{t}$ is a 
$C^{1}$-family of bounded operators, the spectral flow $\Sf(D_{t})_{0\leq t\leq 1}$   
counts the net number of eigenvalues of $D_{t}$ crossing 
the origin as $t$ ranges over $[0,1]$ (see \cite{APS:SARG3}). 
The spectral flow depends  only on the endpoints $D_{0}$ and $D_{1}$ 
and we define  
\begin{equation}
    \Sf(D_{0},D_{1}) = \Sf (D_{t})_{0\leq t\leq 1}.
\end{equation}
Here $D-U^{*}DU= U^{*}[D,U]$ is bounded and one can prove that 
 \begin{equation}
     \Sf(D,U^{*}DU)=   \ind P U P. 
     \label{eq:chern.spectral-flow-index}
 \end{equation}

The \emph{cyclic cohomology groups}  $HC^{*}(\cA)$ of the algebra $\cA$ are obtained from the 
spaces $C^k(\cA)=\{\text{$(k+1)$-linear forms on $\cA$}\}$, $k\in \N$,  by restricting the 
Hochschild 
coboundary, 
  \begin{eqnarray}
      b\psi(a^0, \cdots, a^{k+1}) & = &  \sum (-1)^j \psi(a^0, \cdots, a^j a^{j+1}, \cdots, 
      a^{k+1}) 
       \nonumber  \\
       &  & + \ (-1)^{
      k+1}\psi(a^{k+1}a^0, \cdots, a^k), \qquad  a^j \in \cA, 
  \end{eqnarray}
to cyclic cochains, i.e. those satisfying 
\begin{equation}
    \psi(a^1, \cdots, a^k,a^0) =(-1)^k \psi(a^0, a^{1},\cdots, a^k) \qquad a^j \in \cA.
\end{equation}
It can equivalently be  described as the second filtration of the $(b,B)$-bicomplex of (arbitrary) cochains, where $B : C^m (\cA)
\rightarrow C^{m-1}(\cA)$ is given by  
\begin{gather}
    B=AB_{0}, \qquad (A\phi) (a^0, \cdots, a^{m-1})=\sum (-1)^{(m-1)j} \psi(a^j,\cdots, 
    a^{j-1}),  \\ 
    B_{0}\psi (a^0, \cdots, a^{m-1}) =\psi(1, a^0, \cdots, a^{m-1}), \qquad a^j \in \cA. 
\end{gather}

The \emph{periodic cyclic cohomology} is obtained by taking the inductive limit of the groups 
$HC^k(\cA)$, $k\geq 0$, with respect to the periodicity operator given by the cup product with 
the generator of $HC^{2}(\C)$.  In terms of the $(b,B)$-bicomplex this is the cohomology of the short complex
\begin{equation}
    C^{\ev} (\cA) \stackrel{b+B}{\leftrightarrows} 
    C^{\odd}(\cA), \qquad C^{\ev/\odd} (\cA) = \bigoplus_{k \ \text{even}/\text{odd}} 
    C^{k}(\cA),
\end{equation}
 whose cohomology groups are denoted  $HC^\ev(\cA)$ and $HC^\odd(\cA)$. 

There is a pairing between $HC^\ev(\cA)$ and $K_{0}(\cA)$ such that for any  cocycle $\varphi 
=(\varphi_{2k})$ in $C^\ev(\cA)$ and  for any selfadjoint idempotent $e$ in $M_{q}(\cA)$ we have 
\begin{equation}
    \acou{[\varphi]}{[e]} = \sum_{k \geq 0}(-1)^k \frac{(2k)!}{k!}\varphi_{2k}\# \Tr 
    (e,\cdots,e), 
\end{equation}
where $\varphi_{2k}\# \Tr$ is the$(2k+1)$-linear map on $M_{q}(\cA)=M_{q}(\C)\otimes \cA$ given by
\begin{equation}
  \varphi_{2k}\# \Tr(\mu^{0}\otimes a^0, \cdots, \mu^{2k}\otimes a^{2k}) = 
    \Tr(\mu^{0}\ldots \mu^{2k}) \varphi_{2k}(a^0, \cdots, a^{2k}), 
\end{equation}
for $\mu^j \in M_{q}(\C)$ and $a^j \in \cA$.  

The pairing between $HC^\odd(\cA)$ and $K_{1}(\cA)$ is such that 
 \begin{equation}
    \acou{[\varphi]}{[U]} = \frac1{\sqrt{2i\pi}} \sum_{k \geq 0} 
    (-1)^k k! \varphi_{2k+1}\# \Tr 
    (U^{-1},U,\cdots,U^{-1},U), 
\end{equation}
for any $\varphi =(\varphi_{2k+1})$ in $C^\odd(\cA)$ and  any $U$ in $U_{q}(\cA)$.

\begin{example}%
    Let $\cA$ be the algebra $C^\infty(M)$ of smooth functions on  a compact manifold of 
    dimension $n$ and let  $\cD_{k}(M)$ denote the space of $k$-dimensional de Rham current on $M$. Any 
    $C\in \cD_{k}(M)$ 
    define a Hochschild cochain on $C^\infty(M)$ by letting 
    \begin{equation}
  \psi_{C}(f^{0},f^{1}, \ldots, f^n)= \acou{C}{f^{0}df^{1}\wedge \ldots \wedge df^k} 
        \qquad f^j \in C^\infty(M). 
    \end{equation}
     This cochain satisfies $B\psi_{C}=k\psi_{d^t C}$, where $d^t$ is the de Rham 
    boundary for currents. Thus the map
    \begin{equation}
        \cD_{\ev/\odd}(M)\ni C=(C_{k}) \longrightarrow \varphi_{C}=(\frac1{k!}\psi_{C_{k}}) \in C^{\ev/\odd} (C^\infty(M))
        \label{eq:LIFNCG.morphism-homology-cyclic}
    \end{equation} 
    induces a morphism from the de Rham's homology group $H^{\ev/\odd}(M)$ to the cyclic cohomology group 
    $HC^{\ev/\odd}(C^{\infty}(M))$. This is actually  an 
     isomorphism if we restrict ourselves to the cohomology of continuous cyclic cochains~\cite{Co:NCG}. 
     
     Moreover,  under the Serre-Swan isomorphism 
     $K_{*}(C^\infty(M)) \simeq K^{-*}(M)$ we have, in the even case, 
     \begin{equation}
         \acou{[\varphi_{C}]}{\cE} = \acou{C}{\Ch^{*}_{\ev}\cE} \qquad \forall \cE \in K^{0}(M), 
         \label{eq:chern.acou-even}
     \end{equation}
    where $\Ch^{*}_{\ev}$ is the even Chern character in cohomology 
    (\emph{cf.} Theorem~\ref{thm:aymptotic-AS}),  while in the odd case we have 
    \begin{equation}
            \acou{[\varphi_{C}]}{[U]} = \frac1{\sqrt{2i\pi}} \acou{C}{\Ch^{*}_{\odd}[U]} \qquad 
        \forall U \in C^{\infty}(M,U_{N}(\C)), 
         \label{eq:chern.acou-odd}
    \end{equation}
    where $\Ch^{*}_{\odd}[U] $ is the Chern character of $[U]\in K^{-1}(M)$, i.e. the cohomology class of the odd 
    form $\Ch U=\sum (-1)^{k} \frac{k!}{(2k+1)!} \Tr (U^{-1}dU)^{2k+1}$. 
\end{example}

The index maps~(\ref{eq:LIFNCG.index-even}) and~(\ref{eq:LIFNCG.index-odd}) 
can be computed by pairing $K^{*}(\cA)$ with a cyclic cohomology class 
as follows. Suppose first that the spectral triple $(\cA,\cH,D)$ is \emph{$p$-summable}, i.e. 
 \begin{equation}
     \mu_{k}(D^{-1}) = \op{O}(k^{-1/p})  \qquad \text{as $k\rightarrow +\infty$},
 \end{equation}
 where $ \mu_{k}(D^{-1})$ is the $(k+1)$'th characteristic value of the compact operator 
$D^{-1}$. Then let $\Psi_D^0(\cA)$ denote the algebra generated by the $\delta^k(a)$'s, $a \in 
\cA$, where $\delta$ is the derivation $\delta(T)=[|D|,T]$ (assuming $\cA$  is contained 
in $\cap_{k\geq 0} \op{dom} \delta^k$).

\begin{definition}%
The dimension spectrum of $(\cA,\cH,D)$ is the union set of the singularities of all the zeta 
functions 
$\zeta_{b}(z) = \Tr b|D|^{-z}$, $b\in \Psi_D^0(\cA)$. 
\end{definition}

Assuming simple and discrete dimension spectrum we define an analogue of the 
Wodzicki-Guillemin residue (\cite{Wo:LISA}, \cite{Gu:NPWF}) on $\Psi_D^0(\cA)$ by 
letting 
\begin{equation}
    \bint b = \Res_{z=0} \Tr b|D|^{-z} \qquad \text{for}\ b\in \Psi_D^0(\cA).
\end{equation}
This functional is a trace on the algebra $\Psi_D^0(\cA)$ and is local in the sense of  noncommutative 
geometry since it vanishes on any element of $\Psi_D^0(\cA)$ which is traceable. 

\begin{theorem}[{\cite[Thm.~II.3]{CM:LIFNCG}}]\label{thm:chern.connes-moscovici.even}
Suppose that $(\cA,\cH,D)$  is even, $p$-summable and has a discrete and simple 
dimension spectrum. Then: \smallskip 
 
1) The following formulas define an even cocycle 
$\varphi_{\op{CM}}^\ev =(\varphi_{2k})$ 
        in the $(b,B)$-complex of the algebra $\cA$. For $k=0$, 
             \begin{equation}
    \varphi_{0}(a^{0})  = \text{finite part of $\Tr \gamma a^{0}e^{-tD^{2}}$ as $t\rightarrow 
    0^{+}$},
	\label{eq:chern.connes-moscovici.constant}
    \end{equation}
    while for $k\neq0$, 
    \begin{equation}
        \varphi_{2k}(a^{0}, \ldots, a^{2k}) =   \sum_{\alpha}c_{k,\alpha}\bint  \gamma a^0 [D,a^1]^{[\alpha_{1}]} \ldots 
          [D,a^{2k}]^{[\alpha_{2k}]} |D|^{-2(|\alpha|+k)}, 
          \label{eq:chern.connes-moscovici.even}
       \end{equation}    
    where $\Gamma(|\alpha|+k)c_{k,\alpha}^{-1}=  
    2(-1)^{|\alpha|}\alpha!(\alpha_{1}+1) \cdots (\alpha_{1}+\cdots 
    +\alpha_{2k}+ 2k)$ and the symbol $T^{[j]}$ denotes the $j$'th iterated commutator with 
    $D^2$. \smallskip 

    2)  We have $\ind_{D}(\cE) = \acou{[\varphi_{\op{CM}}]}{\cE}$ for any $\cE\in 
K_{0}(\cA)$. 
\end{theorem}
 
\begin{theorem}[{\cite[Thm.~II.2]{CM:LIFNCG}}]\label{thm:chern.connes-moscovici.odd}
Assume  $(\cA,\cH,D)$ is $p$-summable and has a discrete and simple 
dimension spectrum. Then: \smallskip 

 1) We define an odd  cocycle 
$\varphi_{\op{CM}}^\odd =(\varphi_{2k+1})$ 
        in the $(b,B)$-complex of the algebra $\cA$ by letting 
     \begin{eqnarray}
             \lefteqn{ \varphi_{2k+1}(a^{0}, \ldots, a^{2k+1}) =  } \nonumber \\ 
            &&  \sqrt{2i\pi}
               \sum_{\alpha}c_{k,\alpha}\bint  a^0 [D,a^1]^{[\alpha_{1}]} \ldots 
              [D,a^{2k+1}]^{[\alpha_{2k+1}]} |D|^{-2(|\alpha|+k)-1)},
             \label{eq:chern.connes-moscovici.odd}
       \end{eqnarray} 
        where $ \Gamma(|\alpha|+k +\frac12)c_{k,\alpha}^{-1}= (-1)^{|\alpha|}\alpha!(\alpha_{1}+1) \cdots (\alpha_{1}+\cdots 
    +\alpha_{2k}+2k +1)$. \smallskip	    
 
    2) We have $\ind_{D}(U) = \acou{[\varphi_{\op{CM}}]}{U}$ for any $U\in K_{1}(\cA)$. 
 \end{theorem}

 \begin{example}%
Let $M$ be a compact manifold of dimension $n$ and let $D$ be a  pseudodifferential 
operator of order 1 on $M$  acting on the sections of a vector bundle $\cE$ over $M$ such that $D$ is  
elliptic and selfadjoint. Then the triple 
\begin{equation}
    (C^{\infty}(M),L^{2}(M,\cE),D)
\end{equation}
is an $n$-summable spectral triple, which is  even when $\cE$ is equipped with a 
$\Z_{2}$-grading anticommuting with $D$. 
In any case the  algebra $\Psi_D^0(C^{\infty}(M))$ is contained  in the algebra of 
\pdo's with order $\leq 0$. So by the very construction of the Wodzicki-Guillemin 
residue~(\cite{Wo:LISA},~\cite{Gu:NPWF}) this spectral triple 
has a simple and discrete dimension 
spectrum contained in $ \{k \in \Z; \ k\leq n\}$. 

In fact, given $P\in \Psi^m(M,\cE)$, $m\in\Z$, the function $z \rightarrow 
\Trace P|D|^{-z}$ has a meromorphic continuation on $\C$ with at worst simple 
poles at integers $k$, $k\leq m+n$. At $z=0$ the residue coincides with the Wodzicki-Guillemin
residue $\bint P$ of $P$,   i.e. 
\begin{equation}
    \bint P=\res_{z=0} \Tra P|D|^{-z}=\int_{M}\tr_{\cE} c_{P}(x),
     \label{eq:LIFNCG-NCR}
\end{equation}
where $c_{P}(x)$ is an $\END\cE$-valued density on $M$. 
Hence the  formulas for the CM-cocycle $\varphi_{\op{CM}}$ hold using 
the Wodzicki-Guillemin residue as residual trace. 
 \end{example}
  
\section{The CM cocycle of a Dirac spectral triple (even case)}
\label{sec.even-chern}
 
Let $(M^n,g)$ be a compact Riemannian spin manifold  of even dimension  
and let $\sD_{M}$ denote the Dirac operator acting its spin bundle $\sS$. Then the spectral triple 
$(C^\infty(M), L^{2}(M,\sS),\sD_{M})$ 
is  even and has a discrete and simple dimension spectrum. In this section we shall compute the associated even CM 
cocycle and explain how this allows us to recover the index formula of Atiyah and Singer. 
    
\begin{theorem}\label{thm:even-chern}
    The components of the even 
    CM cyclic cocycle $\varphi_{\op{CM}}^\ev=(\varphi_{2k})$ associated to the spectral triple 
    $(C^\infty(M), L^{2}(M,\sS),\sD_{M})$ are given by  
    \begin{equation}
          \varphi_{2k}(f^{0}, \ldots, f^{2k}) = \frac1{(2k)!}
	  \int_{M} f^{0} df^{1}\wedge \cdots \wedge df^{2k} \wedge 
          \hat{A}(R_{M})^{(n-2k)} , 
	  \label{eq:even-chern.cocycle-Dirac}
    \end{equation}
    for $f^{0},f^{1},\ldots,f^{n}$ in $C^{\infty}(M)$. 
\end{theorem}
\begin{proof}
First, it follows from Theorem~\ref{thm:index.heat-kernel-asymptotic} that  
\begin{equation}
    \varphi_{0}(f^{0})= \lim_{t \rightarrow 0^{+}} \Str f^{0}e^{-t\sD_{M}^{2}} = 
            \int_{M} f^{0}  \hat{A}(R_{M})^{(n)}. 
\end{equation}
Second, let $\alpha$ be a $2k$-fold index, $k\geq 1$, and define 
\begin{equation}
    \sP_{\alpha}=f^{0}[\sD_{M}, f^1]^{[\alpha_{1}]} \cdots [\sD_{M}, f^{2k}]^{[\alpha_{2k}]} 
    = f^{0}c(df^{1})^{[\alpha_{1}]} \cdots c(df^{2k})^{[\alpha_{2k}]}. 
\end{equation}
Then in order to use Formula~(\ref{eq:chern.connes-moscovici.even}) for 
$\varphi_{2k}(f^{0},\ldots,f^{2k})$ we need to compute 
\begin{equation}
    \bint \gamma \sP_{\alpha} |\sD_{M}|^{-2(k+|\alpha|)} =
    \res_{z=0} \Str \sP_{\alpha} 
            |\sD_{M}|^{-2(k+|\alpha|)-z}. 
\end{equation}
The main step is to prove the lemma below.  
\begin{lemma}\label{lem:even-chern.key-lemma}
 For $t>0$ let $k_{\alpha,t}(x,y)$ be the kernel of $\sP_{\alpha} e^{-t\sD_{M}^{2}}$. Then as $t\rightarrow 0^{+}$ we 
 have the following asymptotics   in $C^{\infty}(M,|\Lambda|(M))$: \smallskip 
 
 - $ \Str_{x}  k_{\alpha,t}(x,x)=\op{O}(t^{-(k+|\alpha|)+1})$ if $\alpha\neq 0$;  \smallskip  

 - $ \Str_{x} k_{0,t}(x,x)=  \frac{t^{-k}}{(2i\pi)^{\frac{n}2}} f^{0}df^{1}\wedge \ldots \wedge 
             df^{2k} \wedge \hat{A}(R_{M})^{(n-2k)} + \op{O}(t^{-k+1})$. 
 
 \end{lemma}
\begin{proof}
    In synchronous normal coordinates  $c(df^{j})$ and $\sD^{2}$ have respectively Getzler orders 1 and 2 and model operators 
    $df^{j}(0)$ and  $H_{R}=-\sum 
    (\partial_{i}-R_{ij}^{M}(0)x^j)^{2}$. Therefore, by 
    Lemma~\ref{lem:index.top-total-order-symbol-composition} the operator  $\sP_{\alpha}$ has  
    Getzler order $\leq 2(k+|\alpha|)$ and we have   
    \begin{equation}
       \sP_{\alpha}= c[f^{0}(0)df^{1}(0)^{[\alpha_{1}]}\wedge \cdots \wedge df^{2k}(0)^{[\alpha_{2k}]}] 
        +\op{O}_{G}(2(k+|\alpha|)-1),
    \end{equation}
    where $T^{[j]}$ denotes the $j$'th iterated commutator of $T$ with $H_{R}$.  
    Remark that $[H_{R},df^j(0)]=0$, so if $\alpha \neq 0$ then $ \sP_{\alpha}Q$ has 
    Getzler order $\leq 2(k+|\alpha|)-1$. 
    Moreover as the model operator of $P_{0}$ is 
    $\sP_{0(2k)}=f^{0}(0)df^{1}(0)\wedge \ldots \wedge 
         df^{2k}(0)$ we get $(\sP_{0(2k)}G_{R})(0,1)=f^{0}(0)df^{1}(0)\wedge \ldots \wedge 
         df^{2k}(0)\wedge \hat{A}(R^{M}(0))$. The result then follows by applying 
         Proposition~\ref{prop:AS.differantiated-local-index-theorem}.  
\end{proof}

Now, by the Mellin formula  we have
$\sD_{M}^{-2s}=\Gamma(s)^{-1}\int_{0}^{\infty} t^{s-1} e^{-t\sD^{2}_{M}}dt$ for $\Re s>1$, so
the function $\Str \sP_{\alpha}|\sD_{M}|^{-(2+|\alpha|)-2z}$ 
 coincides with 
 \begin{equation}
   \Gamma(k+|\alpha|+z) ^{-1} \int_{0}^1 t^{k+|\alpha|+z} \Str (\sP_{\alpha} 
 e^{-t\sD^{2}_{M}})  \frac{dt}{t},
 \end{equation}
 up to a holomorphic function 
 on the halfplane $\Re z>-1$. 
Therefore, it follows from 
Lemma~\ref{lem:even-chern.key-lemma} that if $\alpha \neq 0$ then  $\Str 
\sP_{\alpha}|\sD_{M}|^{-(2+|\alpha|)-2z}$ has an analytic  continuation on the halfplane $\Re z>-1$, while $\Str 
\sP_{0}|\sD_{M}|^{-2-2z}$  is equal to   
\begin{equation}
 \frac{(2i\pi)^{-\frac{n}2}}{z\Gamma(z+k)} \int_{M}f^{0}df^{1}\wedge \ldots \wedge 
         df^{2k} \wedge \hat{A}(R_{M})^{(n-2k)}, 
\end{equation}
modulo a holomorphic function on the halfplane $\Re z>-1$. Thus in the formula~(\ref{eq:chern.connes-moscovici.even}) for 
 $\varphi_{2k}(f^{0},\ldots,f^{2k})$ all the 
residues corresponding to 
$\alpha \neq 0$ are zero, while for $\alpha =0$ we get
\begin{equation}
    \bint \gamma \sP_{0} |\sD_{M}|^{-2k} = \frac{2(2i\pi)^{-\frac{n}2}}{(k-1)!}
	  \int_{M} f^{0} df^{1}\wedge \ldots \wedge df^{2k} \wedge 
          \hat{A}(R_{M})^{(n-2k)}. 
\end{equation}
This gives $ \varphi_{2k}(f^{0}, \ldots, f^{2k}) = \frac{(2i\pi)^{-\frac{n}2}}{(2k)!}
	  \int_{M} f^{0} df^{1}\wedge \ldots \wedge df^{2k} \wedge 
          \hat{A}(R_{M})^{(n-2k)}$  since 
$c_{k,0}=\frac12 \frac{\Gamma(k)}{(2k)!}$.  
\end{proof}

We  can now recover the local index formula of Atiyah and Singer. 
Let $\cE$ be a Hermitian vector bundle over $M$ together with a unitary connection with 
curvature $F^\cE$ and let $\sD_{\cE}$ denote the associated twisted Dirac operator. The starting point is that under 
the Serre-Swan isomorphism $K^{0}(M)\simeq K_{0}(C^{\infty}(M))$ we have $\ind \sD_{\cE}^{+}=\ind_{\sD_{M}}[\cE]$ 
(e.g.~\cite[Sect.~2]{Mo:IEDPNG}). Therefore, from Theorem~\ref{thm:chern.connes-moscovici.even} we obtain 
\begin{equation}
    \ind \sD_{\cE}^{+}=\acou{[\varphi_{CM}]}{[\cE]}. 
\end{equation}

On the other hand, Formula~(\ref{eq:even-chern.cocycle-Dirac}) shows that $\varphi_{CM}$ is the image under the 
map~(\ref{eq:LIFNCG.morphism-homology-cyclic}) of the even de Rham current 
that is the Poincar\'e dual of $\hat{A}(R_{M})$. Thus using~(\ref{eq:chern.acou-even}) we get
\begin{equation}
 \ind \sD_{\cE}^{+}= (2i\pi)^{-\frac{n}2} \int_{M}
[\hat{A}(R^{M})\wedge \Ch F^{\cE}]^{(n)},   
\end{equation}
which is precisely the index formula of Atiyah and Singer.  

\section{The CM cocycle of a Dirac spectral triple (odd case)}
\label{sec.odd-chern}
 In this section we 
compute the CM cocycle  corresponding to a 
Dirac operator on an odd dimensional spin manifold. As consequence we  can recapture the spectral 
flow formula of Atiyah-Patodi-Singer~\cite{APS:SARG3}. 

Let $(M^n,g)$ be a compact Riemannian spin manifold of odd dimension and let $\sS$ be a spin 
bundle for $M$, so that each fiber $\sS_{x}$ is an irreducible representation space for 
$\Cl_{x}(M)$. The Dirac operator $\sD_{M}$ acting on the sections of $\sS$ is given by the 
composition,
\begin{equation}
    C^\infty(M, \sS) \stackrel{\nabla^{\sS}}{\longrightarrow} C^\infty(M,T^{*}M \otimes\sS) 
    \stackrel{c\otimes 1}{\longrightarrow} C^\infty(M,\sS),
     \label{eq:odd-chern.Dirac-operator}
\end{equation}
where  $c$ denotes the action of $\Lambda T^{*}M$ on $\sS$ by Clifford 
representation.  This gives rise to an odd spectral triple $(C^\infty(M), L^{2}(M,\sS),\sD_{M})$ with simple and 
discrete dimension spectrum. 

\begin{theorem}\label{thm:odd-chern}
The components of the odd 
    CM cocycle $\varphi_{\op{CM}}^\odd=(\varphi_{2k+1})$ associated to the spectral triple 
    $(C^\infty(M), L^{2}(M,\sS),\sD_{M})$ are given by 
    \begin{eqnarray}
        \lefteqn{\varphi_{2k+1}(f^{0}, \ldots, f^{2k+1}) = } \nonumber  \\ 
&& 	  \sqrt{2i\pi} 
\frac{(2i\pi)^{-[\frac{n}2]+1}}{(2k+1)!} 
	  \int_{M} f^{0} df^{1}\wedge \ldots \wedge df^{2k+1} \wedge 
          \hat{A}(R_{M})^{(n-2k-1)}.
	  \label{eq:odd-chern.cocycle-Dirac}
    \end{eqnarray}
     for $f^{0},\ldots,f^{(n)}$ in $C^\infty(M)$. 
\end{theorem}

Before tackling the proof of Theorem~\ref{thm:odd-chern} let us explain the similarities with the even case. 
Since the dimension of $M$ is odd there is not anymore an isomorphism between 
$\Cl(M)$ and $\End \sS$ and so we 
need to distinguish between them. In fact if $e_{1}, \ldots, e_{n}$ is an orthonormal frame for 
$T_{x}M$ then 
$c(e^{1})\cdots c(e^{n})$ acts like $(-i)^{[\frac{n}2]+1}$ on $\sS_{x}$ (\emph{cf.}~\cite{Ge:POSASIT}, \cite{BF:AEF2}). 
Therefore, if we look at 
$c(e^{i_{1}})\cdots c(e^{i_{k}})$, $i_{1}<\ldots<i_{k}$, as an endomorphism of $\sS_{x}$ then we 
have 
\begin{equation}
    \tr_{\sS_{x}} c(e^{i_{1}})\cdots c(e^{i_{k}}) = \left\{  \begin{array}{ll}
                                                              0 & \text{if $0<k< n$,}   \\
          (-i)^{[\frac{n}2]+1} 2^{[\frac{n}2]}& \text{if $k=n$}. 
    \end{array} \right. 
    \label{eq:odd-chern.trace-clifford}
\end{equation}
Therefore, provided only an odd number of Clifford variables are involved 
the trace behaves as the supertrace in even dimension. 

Bearing this in mind let $Q\in \Psi_{\op{v}}^{*}(M\times \R, \sS)$. In synchronous normal coordinates $Q$ is given near 
the origin by a Volterra \psido\ operator $\tilde{Q}$ which acts on the trivial bundle with fiber $\sS_{n}$, i.e.~is 
with coefficients in $\End \sS_{n}$. Thus, via the Clifford representation  $\tilde{Q}$ comes from 
a Volterra \psido\ ${}^\Cl Q$ on $\R^{n}\times \R$ with coefficients in $\Cl^\odd(\R^n)$. Then 
using~(\ref{eq:odd-chern.trace-clifford}) we get 
   \begin{equation}
       \Tr_{} K_{Q}(0,0,t) = -i(-2i)^{[\frac{n}2]} 
       \sigma[ K_{{}^{\Cl}\sP{}^\Cl Q_{1}}(0,0,t)]^{(n)}+ \op{O}(t^{\infty}).  
       \label{Odd-CM.trace-Kalpha}
   \end{equation}
   
 On the other hand, in the the proof of Theorem~\ref{thm:index.heat-kernel-asymptotic} we identified $\End 
  \sS_{n}$ and $\Cl(\R^{n})$. Thus definition~\ref{def:getzler.model-operator} and  all the 
  lemmas~\ref{lem:AS.approximation-asymptotic-kernel}-\ref{lem:getzler-model-asymptotic}  
  hold \emph{verbatim} for 
  Volterra \psido's with coefficients in $\Cl(\R^{n})$, and this independently of the parity of $n$. For instance, if we 
  let $m$ be the Getzler order of ${}^{\Cl}Q$ then from Lemma~\ref{lem:AS.approximation-asymptotic-kernel} we get: \smallskip 
  
  - $\sigma[K_{{}^{\Cl}Q}(0,0,t)]^{(j)}= \op{O}(t^{\frac{j-m-n-1}2})$ if $m-j$ is odd, \smallskip 
  
  - $\sigma[ K_{{}^{\Cl}Q}(0,0,t)]^{(j)}=t^{\frac{j-m-n}2-1} 
         K_{{}^{\Cl}Q_{(m)}}(0,0,1)^{(j)} + \op{O}(t^{\frac{j-m-n}2}) $ otherwise.\smallskip 
  
  \noindent Note this is consistent with~(\ref{eq:AS-convergence-symbol-KQ}) because as $n$ is odd by the proof of 
  Lemma~\ref{lem:AS.approximation-asymptotic-kernel} whenever $m$ is even we have 
\begin{equation}
    K_{{}^{\Cl}Q_{(m)}}(0,0,1) = \sigma[\check{q}_{-2-n}(0,0,1)]^{(n)}=0 .  
    \label{eq:CM-odd.vanishing-K-model}
\end{equation}
\begin{definition}%
    We say that $Q\in \Psi_{\op{v}}^{*}(M\times \R, \sS\otimes \cE)$ has Getzler order $m$ if in synchronous normal coordinates 
    centered at any $x_{0}\in M$ the operator ${}^{\Cl}Q$ defined as above has Getzler order $m$. Moreover we let 
    $Q_{(m)}={}^{\Cl}Q_{(m)}$ be the model operator of $Q$.
\end{definition}

Along similar lines as that of the proof of Proposition~\ref{prop:AS.differantiated-local-index-theorem} we obtain: 

\begin{proposition}\label{prop:AS.differantiated-local-index-theorem-odd}
    Let $\sP$ be a differential operator on $M$ acting on $\sS$ with Getzler order $m$ and 
    let $h_{t}(x,y)$ denote the kernel of $\sP e^{-t\sD_{\cE}^{2}}$. Then as $t\rightarrow 0^{+}$ we 
    have an asymptotics in $C^{\infty}(M,  |\Lambda|(M))$ of the form: \smallskip 
 
        - $\Tr_{x}h_{t}(x,x)= \op{O}(t^{\frac{-m+1}2})$ if $m$ is even;  \smallskip  
  
        -  $\Tr_{x}h_{t}(x,x)= t^{\frac{-m}2} 
        B_{0}(\sD^{2}_{\cE},\sP)(x) + \op{O}(t^{\frac{-m}2+1}) $ if  $m$ is odd, where  in synchronous normal coordinates 
        centered at $x_{0}$ and with $\sP_{(m)}$ denoting the model operator of $\sP$  we have 
        $B_{0}(\sD^{2}_{\cE},\sP)(0)= (-i)^{[\frac{n}2]+1} 2^{[\frac{n}2]} [(\sP_{(m)}G_{R})(0,1)]^{(n)}$. 
\end{proposition}

\begin{proof}[Proof of Theorem~\ref{thm:odd-chern}]
Let  $\sP_{\alpha}= 
    f^{0}[\sD_{M}, f^1]^{[\alpha_{1}]} \cdots [\sD_{M}, f^{2k+1}]^{[\alpha_{2k+1}]}$ where $\alpha$ is a 
    $(2k+1)$-fold index. Then applying Proposition~\ref{prop:AS.differantiated-local-index-theorem-odd} 
    and arguing as in the proof of Lemma~\ref{lem:even-chern.key-lemma} 
    shows that 
as $t \rightarrow 0^{+}$ we have: \smallskip 

- $\Tr \sP_{\alpha}e^{-t\sD_{M}^{2}}= \op{O}(t^{-(k+|\alpha|)+\frac12})$ if $\alpha\neq 0$,\smallskip 

- $\Tr \sP_{0}e^{-t\sD_{M}^{2}}= 
\frac{t^{-k-\frac12}}{(2i\pi)^{[\frac{n}2]}} \int_{M} f^{0} df^{1}\wedge \ldots \wedge df^{2k+1} \wedge 
          \hat{A}(R_{M})^{(n-2k-1)}+ \op{O}(t^{-k+\frac12})$.\smallskip 
          
\noindent Then as in the proof of 
Theorem~\ref{thm:even-chern} we deduce that in the formula~(\ref{eq:chern.connes-moscovici.odd}) for 
$\varphi_{2k+1}(f^{0}, \ldots , f^{2k+1})$ only  $\bint \sP_{0} |\sD_{M}|^{-(2k+1)}$ is  nonzero and equal to 
\begin{equation}
        \frac2{\Gamma(k+\frac12)} \frac{(2i\pi)^{-[\frac{n}2]}}{2i\sqrt{\pi}}     
        \int_{M} f^{0} df^{1}\wedge \ldots  \wedge df^{2k+1} \wedge 
              \hat{A}(R_{M})^{(n-2k-1)}.
\end{equation}
Hence $ \varphi_{2k+1}(f^{0}, \ldots , f^{2k+1})=\sqrt{2i\pi} \frac{(2i\pi)^{-[\frac{n}2]+1}}{(2k+1)!}\int_{M} f^{0} 
df^{1}\wedge \ldots \wedge df^{2k+1} 
\wedge  \hat{A}(R_{M})^{(n-2k-1)}$.
 \end{proof}

As a consequence of Theorem~\ref{thm:odd-chern} we can recover the spectral flow formula of 
Atiyah-Patodi-Singer~\cite{APS:SARG3} in the 
case of a Dirac operator (see also~\cite{Ge:OCCCHSF}). 
\begin{theorem}[{\cite[p. 95]{APS:SARG3}}]\label{thm:odd-chern.APS}
 For any $U \in C^\infty(M, U(N))$ we have 
 \begin{equation}
     \Sf(\sD_{M}, U^{*}\sD_{M} U)= (2i\pi)^{-[\frac{n}2]-1} 
     \int_{M}[ \hat A (R^M)\wedge \Ch (U)]^{(n)}.
      \label{eq:odd-chern.APS}
 \end{equation}
\end{theorem}
\begin{proof}
    Thanks to~(\ref{eq:chern.spectral-flow-index})  and Theorem~\ref{thm:chern.connes-moscovici.odd} we have 
\begin{equation}
        \Sf(\sD_{M}, U^{*}\sD_{M} U) = \ind_{\sD_{M}}[U]=\acou{[\varphi_{CM}]}{[U]} .
\end{equation}
   Moreover, Formula~(\ref{eq:odd-chern.cocycle-Dirac})
   shows that $\varphi_{CM}$ is the image under the map~(\ref{eq:LIFNCG.morphism-homology-cyclic}) 
   of the odd de Rham current 
that is the Poincar\'e dual of $\hat{A}(R_{M})$. Formula~(\ref{eq:odd-chern.APS}) then follows by using~(\ref{eq:chern.acou-odd}).
\end{proof}

\begin{acknowledgements} 
  I would like to thank Daniel Grieser, Thomas Krainer, Xiaonan Ma, Richard Melrose and Henri Moscovici for helpful and 
  stimulating discussions. 
  This work was partially supported by the nodes of the European Research Training Network 
  HPCRN-CT-1999-00118 "Geometric Analysis'' at Postdam University and at Humboldt University at Berlin and by the 
  NSF collaborative grant DMS 9800765 of N. Higson and J. Roe, when I was visiting at Penn State University.   
\end{acknowledgements}


\begin{thebibliography}{CFKS}
 \bibitem[ABP]{ABP:OHEIT} Atiyah, M., Bott, R., Patodi, V.: 
\emph{On the heat equation and the index theorem.} 
Invent. Math. \textbf{19},  279--330 (1973). 
     
\bibitem[AS1]{AS:IEO1}  Atiyah, M.,  Singer, I.: 
\emph{The index of elliptic operators. I. }
Ann. of Math. (2) \textbf{87},    484--530 (1968).  

\bibitem[AS2]{AS:IEO3}  Atiyah, M.,  Singer, I.: 
\emph{The index of elliptic operators.  III.}
Ann. of Math. (2) \textbf{87},    546--604 (1968).  


\bibitem[APS]{APS:SARG3}  Atiyah, M., Patodi, V., Singer, I.:
\emph{Spectral asymmetry and Riemannian geometry. III. }
Math. Proc. Camb. Philos. Soc. \textbf{79},   71--99 (1976). 

\bibitem[BGS]{BGS:HECRM} Beals, R., Greiner, P.,   Stanton, N.:   
\emph{The heat equation on a CR manifold.} J. Differential Geom. \textbf{20}, 343--387 (1984).

\bibitem[BGV]{BGV:HKDO} Berline, N., Getzler, E., Vergne, M.: 
\emph{Heat kernels and Dirac operators}. 
Springer-Verlag, Berlin, 1992. 

 \bibitem[Bi]{Bi:ASITPA} Bismut, J.-M.: 
\emph{The Atiyah-Singer theorems: a probabilistic approach. I. The index theorem.}
 J. Funct. Anal. \textbf{57}, 56--99 (1984) . 

  \bibitem[BF]{BF:AEF2} Bismut, J.-M., Freed, D.: \emph{The analysis of elliptic families.}
  Comm. Math. Phys \textbf{107},  103--163 (1986).
 
\bibitem[CH]{CH:ECCIDO} Chern, S., Hu, X.:  \emph{Equivariant Chern character for the invariant 
Dirac operator}. Michigan Math J. \textbf{44},  451--473 (1997). 
 
\bibitem[Co]{Co:NCG} Connes, A.:  \emph{Noncommutative geometry}.  
Academic Press,  San Diego, 1994.

\bibitem[CM]{CM:LIFNCG}  Connes, A., Moscovici, H.: \emph{The local index formula in
noncommutative geometry}. GAFA  \textbf{5},  174--243 (1995). 

\bibitem[CFKS]{CFKS:SOAQM} Cycon, H. L.; Froese, R. G.; Kirsch, W.; Simon, B.: 
\emph{Schr\"odinger operators with application to quantum mechanics and global geometry}. 
Texts and Monographs in Physics. Springer-Verlag, Berlin, 1987. 

\bibitem[Ge1]{Ge:POSASIT} Getzler, E.: \emph{Pseudodifferential operators on supermanifolds and 
the Atiyah-Singer index theorem}. Comm. Math. Phys. \textbf{92},  163--178 (1983). 

\bibitem[Ge2]{Ge:SPLASIT} Getzler, E.: 
\emph{A short proof of the local Atiyah-Singer index theorem}. 
Topology \textbf{25}, 111--117 (1986). 

\bibitem[Ge3]{Ge:OCCCHSF} Getzler, E.: The odd Chern character in cyclic homology and 
spectral flow. Topology \textbf{32},  489--507 (1993).  

\bibitem[Gi]{Gi:ITHEASIT} Gilkey, P.: 
\emph{Invariance theory, the heat equation, and the Atiyah-Singer index theorem}. 
 Publish or Perish,  1984.

\bibitem[GJ]{GJ:QPFIPV} Glimm, J.; Jaffe, A.: \emph{Quantum physics. A functional integral point of view.}
Springer-Verlag,  1987.

\bibitem[Gr]{Gr:AEHE} Greiner, P.: \emph{An asymptotic expansion for the heat equation}. Arch. Rational Mech. Anal. 
\textbf{41},  163--218 (1971).    

\bibitem[Gu]{Gu:NPWF} Guillemin, V.:
\emph{A new proof of Weyl's formula on the asymptotic distribution of eigenvalues}. 
Adv. in Math. \textbf{55},  131--160 (1985). 

\bibitem[Ha]{Ha:LCPLPDE} ÊHadamard, J.: \emph{Lectures on Cauchy's problem in linear partial differential equations.} 
Dover Publications, 1953. 

\bibitem[Hi]{Hi:LIFNCG} Higson, N.: \emph{The local index formula in noncommutative geometry.} Lectures given at the CIME 
Summer School and Conference on algebraic $K$-theory and its applications, Trieste, 2002. Preprint, 2002.

\bibitem[Le]{Le:TSPSCI} Lescure, J.-M.: \emph{Triplets spectraux pour les vari\'et\'es \`a 
singularit\'e conique isol\'ee}. Bull. Soc. Math. France \textbf{129},  593--623 (2001). 

\bibitem[LM]{LM:SG} Lawson, B., Michelson, M.-L.: \emph{Spin Geometry}. Princeton Univ. Press, Princeton,
1993. 

\bibitem[Me]{Me:APSIT} Melrose, R.: \emph{The Atiyah-Patodi-Singer index theorem}. A.K. Peters, Boston, 1993.  

 \bibitem[Mo]{Mo:IEDPNG} Moscovici, H.: \emph{Eigenvalue inequalities and Poincar\'e duality in noncommutative geometry}. 
Comm. Math. Phys. \textbf{184},  619--628 (1997).

 \bibitem[Pi]{Pi:COPDTV} Piriou, A.: \emph{Une classe d'op\'erateurs pseudo-diff\'erentiels du type de Volterra}. Ann. Inst. 
 Fourier  \textbf{20}, 77--94 (1970).

 \bibitem[Ro]{Roe:EOTAM} Roe, J.:  \emph{Elliptic operators, topology and asymptotic methods}. 
  Pitman Research Notes in Mathematics Series 395,  Longman, 1998. 
 
  \bibitem[Ta]{Ta:PDE2} Taylor, M. E.: \emph{Partial differential equations. II. Qualitative studies of linear equations.}
  Applied Mathematical Sciences, 116. Springer-Verlag, New York, 1996.
  
 \bibitem[Wo]{Wo:LISA} Wodzicki, M.: \emph{Local invariants of spectral asymmetry}. 
Invent. Math. \textbf{75},  143--177 (1984). 
\end{thebibliography}
\end{document}